\documentclass{article}

\usepackage{amssymb}
\usepackage{latexsym}

\usepackage{url}
\usepackage{xcolor}

\usepackage{amsmath,amssymb,mathtools}
\usepackage{enumitem} 
\usepackage{xcolor,ulem,soul} 
\usepackage{graphicx, graphics, subcaption}
\usepackage{hyperref}
\usepackage{todonotes}
\usepackage{algorithmic,algorithm}
\usepackage{booktabs,siunitx} 
\usepackage{hyperref}
\usepackage[letterpaper,top=2cm,bottom=2cm,left=3cm,right=3cm,marginparwidth=1.75cm]{geometry}

\newcommand{\wt}{\widetilde}
\newcommand{\eps}{\varepsilon}

\newcommand{\R}{\mathbb{R}}

\newcommand{\abs}[1]{\left\lvert {#1} \right\rvert}
\newcommand{\norm}[1]{\left\lVert {#1} \right\rVert}

\newcommand{\D}{\Delta}

\newcommand{\cC}{\mathcal{C}}

\newcommand\palla[2][m]{\ifx i#1 \todo[linecolor=violet,backgroundcolor=violet!25,bordercolor=violet,inline]{#2}\else
\todo[linecolor=violet,backgroundcolor=violet!25,bordercolor=violet]{#2}\fi}
\newcommand\pes[2][m]{\ifx i#1 \todo[linecolor=magenga,backgroundcolor=magenta!25,bordercolor=magenta,inline]{#2}\else
\todo[linecolor=magenta,backgroundcolor=magenta!25,bordercolor=magenta]{#2}\fi}
\newcommand\ale[2][m]{\ifx i#1 \todo[linecolor=red,backgroundcolor=red!25,bordercolor=red,inline]{#2}\else
\todo[linecolor=red,backgroundcolor=red!25,bordercolor=red]{#2}\fi}
\newcommand\agn[2][m]{\ifx i#1 \todo[linecolor=cyan,backgroundcolor=cyan!25,bordercolor=cyan,inline]{#2}\else
\todo[linecolor=cyan,backgroundcolor=cyan!25,bordercolor=cyan]{#2}\fi}

\newtheorem{remark}{Remark}

\title{Online identification and control of PDEs\\ via Reinforcement Learning methods}%

\author{Alessandro Alla\footnote{Department of molecular sciences and nanosystems, Università Ca' Foscari Venezia, Italy.
 (alessandro.alla@unive.it)}, Agnese Pacifico\footnote{Department of Mathematics, Sapienza University of Rome, Italy. (agnese.pacifico@uniroma1.it)}, Michele Palladino\footnote{Department of Information Engineering, Computer Science and Mathematics, University of L'Aquila, Italy
 (michele.palladino@univaq.it)}, Andrea Pesare\footnote{Bending Spoons, Milan, Italy (andreapesare1@gmail.com}}

\begin{document}
\maketitle

\begin{abstract}
We focus on the control of unknown Partial Differential Equations (PDEs). The system dynamics is unknown, but we assume we are able to observe its evolution for a given control input, as typical in a Reinforcement Learning framework. We propose an algorithm based on the idea to control and identify on the fly the unknown system configuration. In this work, the control is based on the State-Dependent Riccati approach, whereas the identification of the model on Bayesian linear regression. At each iteration, based on the observed data, we obtain an estimate of the \textit{a-priori} unknown parameter configuration of the PDE  and then we compute the control of the correspondent model. We show by numerical evidence the convergence of the method for infinite horizon control problems.

\end{abstract}





 
\section{Introduction}

Reinforcement Learning (RL) is one of the main Machine Learning
paradigms, together with supervised and unsupervised Learning. In
RL, an agent interacts with an unknown environment,
aiming at an action-selection strategy to optimize a system's performance. Generally speaking, one can consider two main RL philosophies.
The first one, called model-based, usually concerns the reconstruction of a model from the data trying to mimic the unknown environment. That model is then used to plan
and compute a suboptimal policy. The second RL philosophy, called model-free,
employs a direct approximation of the value function and/or a policy based
on a dynamic-programming-like algorithm without using a model to simulate the
unknown environment. Model-free methods include the famous Monte Carlo methods
\cite{SuttonBarto1edn}, Temporal-Difference Learning \cite{sutton1988learning,rummery1994line} and Q-Learning \cite{watkins1989learning} and more recent ones \cite{Mnih2015,TRPO2015,Lillicrap2016,SAC2018}. An overview of the two RL approaches can be found in \cite{SuttonBarto1edn}.



Since the philosophy of this work is to connect RL problems and optimal control, we will first recall some classical approaches to optimal control problems. Specifically, we are interested in feedback control to obtain a state dependent optimal
control which is a valuable property as it makes the control system stable with
respect to random disturbances. 
Dynamic Programming (DP, \cite{Bardi-Capuzzo,falcone2013semi}) considers a family of optimal control problems with different initial conditions and states and looks at the relationship between
these problems. The main ingredient is the value function, defined as the
minimum of the cost functional which is the solution of the Bellman equation: a nonlinear
partial differential equation (PDE) of the Hamilton-Jacobi type. Once the value function has
been obtained, it provides optimal feedback control. Although the DP approach is preferable from a theoretical
point of view as it provides sufficient conditions and synthesis of optimal feedback
control, it has always been challenging to apply it to real problems since
it is very expensive from a computational point of view. It suffers from the so-called
“curse of dimensionality” (an expression coined by Bellman himself \cite{bellman1954}),
which means that the computational cost necessary to solve the Bellman equation
grows exponentially with respect to the system dimension. This led to the development
of suboptimal solution methods to the Bellman equation, known as approximate
dynamic programming \cite{bertsekas2008approximate}, that could mitigate the curse of dimensionality.
Later, we started calling these methods Reinforcement Learning \cite{bertsekas2019reinforcement,SuttonBarto}.
Optimal control and RL are strongly connected, as they deal with
similar problems; in fact, both can be regarded as sequential decision problems,
in which one has to make decisions in sequence, trying to optimize not only the
immediate rewards but also the future, delayed ones. Recently, in \cite{powell2007approximate,powell2021reinforcement}, it has been proposed
a unified framework for all the families of sequential decision problems including
OC and RL. More precisely, RL deals with control problems in which
the system's dynamics is uncertain. \\

In this paper, we want to control an unknown nonlinear dynamics following a RL strategy. We will adopt an online strategy as explained below. We suppose that the system is described by a parametric PDE, whose parameters are unknown. We also assume to have a library which includes all possible terms of the PDE, so that a function of the library enters into the model if the corresponding parameter is not zero. Those parameters are the one we need to discover to achieve our goal. The chosen library, in this work, will contain several models which are very well studied in the mathematical physics community. Furthermore, although the system dynamics is unknown, we assume it is always possible to observe the true evolution of the system for a given control input. The possibility to observe the unknown system is a typical assumption in Reinforcement Learning where an agent takes action based on his observation. This will allow us to update the parameter estimate.

To achieve our goal we propose the following workflow: "control--observe--estimate". 
To set the method into perspectives, we begin with an initial parameter estimate that allows to compute the control for such configuration. Note that the control is computed for a problem that uses a parameter estimate and might be far from the true optimal control. Then, by applying that control, we observe the true system configuration by its trajectories. Thus, to update the parameter estimate, we set a linear system based on the observed trajectories which will be solved via Bayesian Linear Regression methods (see e.g. \cite{Rasmussen2006,rossi2012bayesian}).
We iterate this procedure till the end of the chosen time horizon. We will also discuss a heuristic stopping criteria for the paramter estimation. As mentioned, this is an online approach since we update the parameter estimate every iteration. A first approach driven by the same workflow proposed here has been introduced in \cite{pacifico2021LSSC}
for linear low-dimensional problems and quadratic cost functionals. Here, we extend to generic nonlinear control problems with a keen focus on the control of PDEs. The dimension of the discretized problem increases also the challenges of the problem. Our approach to the control of the PDE is based on the discretization by finite differences that reduces the problem to a large system of ordinary differential equations. In the paper, we also show numerically how the computed control stabilizes the PDE for smaller spatial discretization leading to the control of the continuous PDE.





 

Let us now comment on how we solve the control problems. As already mentioned at the beginning of this section, control in feedback form is usually obtained by the solution of dynamic programming equations \cite{Bardi-Capuzzo} or by Nonlinear Model Predictive Control (NMPC, \cite{larsnmpc}). An alternative, which combines elements from both
dynamic programming and NMPC, is the State-Dependent Riccati Equation
(SDRE) approach (see e.g. \cite{C97,BLT07}). The SDRE method originates from the dynamic
programming associated to infinite horizon optimal stabilization. It circumvents its solution by reformulating the feedback synthesis as the sequential solution of state-dependent Algebraic Riccati Equations (ARE), which are updated online along a trajectory. The SDRE feedback is implemented similarly as
in NMPC, but the online solution of an optimization problem is replaced by a nonlinear matrix equation. Later, in \cite{BH18} it has been studied the method theoretically proving conditions for the stabilitazion of the problem whereas in \cite{AKS23} it has been shown that SDRE is also efficient for large scale problems although the high number of ARE one should~solve.

For the sake of completeness we also recall that system identification of nonlinear dynamics is a very active and modern research area with a vast literature. Although our identification is strictly linked to a control problem, we briefly recall some literature. Clearly, Physics Informed Neural Networks (PINNs) deserves to be mentioned due to its innovative, accurate and efficient way to discover partial differential equations using neural network and using information from system in the definition of the loss function. We refer to e.g. \cite{RPK19,KKLPWY21} for a complete description of the method. It is worth to mention also methods based on variants of sparse optimization techniques such as SINDy for ODEs \cite{BPK16} and for PDEs \cite{RBPK16,RABK19}. SINDy was also applied for the identification of controlled problems (see e.g. \cite{KKB19}). The authors used an external source as input to identify the system and then apply NMPC to control the identified model. There, the authors used the workflow: identify first, control later which is different from the strategy presented in the current work. Other strategies dedicated on control and system identification can be found in e.g. \cite{KS08} for PDEs and in e.g. \cite{MLG20} for ODEs. Recently, a study on the control of unknown problem with MPC has been introduced in \cite{CFDKMRV21}. There, the system was identified using the Extended Dynamic Mode Decomposition, i.e. a surrogate linear model in contrast to our work where we directly identify the nonlinear model.

The outline of the paper is the following. In Section \ref{sec:BLR} we recall the basics of Baysian Linear Regression as a building block when adapting our parameter estimate. In Section \ref{sec: sdre}, we briefly explain the State-Dependent Riccati equation. In Section \ref{sec:idcon}, we provide all the details of the method proposed in this paper. Later, numerical experiments to support our algorithm are presented. Finally, conclusions are driven in Section \ref{sec:conclusions}.

\section{Bayesian linear regression}\label{sec:BLR}
Bayesian Linear Regression (BLR, \cite{Rasmussen2006,rossi2012bayesian}),  is a probabilistic method for solving the classical linear regression (LR, \cite{freedman2009statistical}) problem. In LR, we consider \textit{data} in the form of input-output pairs $$\mathcal{D} = \{(x_i,y_i)\}_{i=1,\dots,d}$$ and we suppose that the output variable $y_i \in \mathbb{R}$ can be expressed approximately as a linear function of the input variable $x_i \in \mathbb{R}^n$, i.e. 
\begin{equation}\label{eq: LR}
    y_i \approx x_i^T \theta, \qquad \textrm{ for } \ i=1, \dots, d .
\end{equation}
We look for a parameter $\theta \in \mathbb{R}^n$ such that \eqref{eq: LR} is satisfied. The \textit{(ordinary) least squares (LS)} approach choose $\theta$ by minimizing the sum of squared residuals 
\begin{equation}\label{eq:SSR}
    E(\theta) = \sum_{i=1}^d \abs{y_i - x_i^T\theta}^2 .
\end{equation}
The \textit{LS solution} can be computed analytically and is given by 
\begin{equation}\label{eq: LS}
    \theta_{LS} = (X^T X)^{-1} X^T Y ,
\end{equation}
where we collected all the observed inputs in a matrix $X \in \mathbb{R}^{d \times n}$ and all the observed outputs in a vector $Y \in \mathbb{R}^d$:

\begin{equation}\label{eq: BLR matrices}
    X = \left(\rule{0cm}{1cm} \begin{array}{cccc}
        \qquad &   x_1^T &      \qquad    \\
        \quad &   x_2^T &      \quad    \\
        \quad &    \vdots &      \quad    \\
         \quad &   x_d^T &      \quad    \\
    \end{array} \right), \quad
    Y = \begin{pmatrix}
        y_1 \\
        y_2 \\
        \vdots \\
        y_d
    \end{pmatrix} .
\end{equation}
In BLR, instead, we assume that the deviation of the data from the linear model can be described by a Gaussian noise $\eps_i \sim \mathcal{N}(0,\sigma^2)$:
\begin{equation}\label{eq: likelihood relation}
    y_i = x_i^T\theta + \eps_i ,    
\end{equation}
where $\theta \in \R^n$ is an unknown parameter to be determined and $\sigma>0$. We will assume that the value of $\sigma$ is known, though more general formulations apply Bayesian inference on $\sigma$ as well. Equation~\eqref{eq: likelihood relation} corresponds to fix a conditional distribution of the random variable $y$ given the value of $x$ and $\theta$,
\begin{equation}\label{eq: single likelihood}
   p(y | x, \theta) \sim \mathcal{N}( x^T\theta, \sigma^2) .
\end{equation}
This is what in Bayesian inference is called the \textit{likelihood function}. If we assume that the $d$ observations are independent, the global likelihood function can be written as
\begin{equation}\label{eq: likelihood}
    p(Y | X, \theta) = \prod_{i=1}^d p(y_i | x_i, \theta) \sim \mathcal{N}(X \theta, \sigma^2 I_d),
\end{equation}
where $X, Y$ have been defined in \eqref{eq: BLR matrices}, and $I_d$ denotes the $d$-dimensional identity matrix.

The available information on the parameter $\theta$ is included in the model through the definition of a \textit{prior distribution}, which we assume to be Gaussian with initial mean 
$m_0 \in \mathbb{R}^n$ and covariance matrix $\Sigma_0 \in \mathbb{R}^{n \times n}$:
\begin{equation}\label{dist:theta}
 \theta \sim \mathcal{N}(m_0,\Sigma_0) . 
\end{equation}
Bayesian formulas allow to compute the \textit{posterior distribution} of the parameter $\theta$, which is again a Gaussian distribution \cite{rossi2012bayesian,box2011bayesian}
\begin{equation}\label{eq: posterior}
    p(\theta | X, Y ) = \frac{p(\theta) p(Y | X, \theta)}{\int_{\mathbb{R}^n} p(\theta') p(Y | X, \theta') d\theta'} \sim \mathcal{N}(m,\Sigma),
\end{equation}
where
\begin{equation}\label{eq: posterior matrices}
    \Sigma^{-1} = \frac{1}{\sigma^2} X^T X + \Sigma_0^{-1} \ \textrm{ and } \ m = \Sigma \left( \frac{1}{\sigma^2} X^T Y + \Sigma_0^{-1} m_0 \right) .
\end{equation}
From the posterior distribution one can extract a point estimate of the parameter $\theta$, that is the posterior mean 
\begin{equation}\label{thetaBLR}
 \bar{\theta}_{BLR} = \Sigma \left( \frac{1}{\sigma^2} X^T Y + \Sigma_0^{-1} m_0 \right) = \left( \frac{1}{\sigma^2} X^T X + \Sigma_0^{-1} \right)^{-1} \left( \frac{1}{\sigma^2} X^T Y + \Sigma_0^{-1} m_0 \right) . 
 \end{equation}
However, the advantage of BLR is that it provides a quantification of the uncertainty of this estimate. Finally, we remark that the estimate $\bar{\theta}_{BLR}$ converges to the LS solution \eqref{eq: LS}, when the noise variance $\sigma$ goes to 0. 

\section{Control of nonlinear problem via State Dependent Riccati Equation}\label{sec: sdre}

In this section, we recall one possible approach to control nonlinear differential equations.
 We consider the following infinite horizon optimal control problem:

\begin{equation}\label{disc_cost}
\underset{u(\cdot)\in \mathcal{U}}{\min}\;J(u(\cdot)):=\int\limits_0^\infty \Big(\norm{x(t)}_Q^2+\norm{u(t)}_R^2\Big)\, dt
\end{equation}
subject to the nonlinear dynamical constraint

\begin{align}\label{eq: dynamics}
\begin{aligned}
	\dot x(t)&=A(x(t))x(t)+B(x(t))u(t),\quad t\in(0,\infty),\\ x(0)&=x_0,
	\end{aligned}
\end{align}
where $x(t):[0,\infty]\rightarrow\R^d$ denotes the state of the system, $A(x):\R^d\rightarrow\R^{d\times d}$, the control signal $u(\cdot)$ belongs to $\mathcal{U}:=L^{\infty}(\R^+;\R^m)$ and $B(x):\R^d\rightarrow\R^{d\times m}$. The running cost is given by $\norm{x}_Q^2:=x^{\top}Qx$ with $Q\in\R^{d\times d},\,Q\succ0$, and $\norm{u}_R^2=u^{\top}Ru$ with $R\in\R^{m\times m},\,R\succ0$. This formulation corresponds to the asymptotic stabilization of nonlinear dynamics towards the origin.

We can synthesize a suboptimal feedback law by following e.g. an approach known as the State-dependent Riccati Equation (SDRE). We refer to e.g. \cite{C97,BLT07} for more details on the topic.

The SDRE approach is based on the idea that infinite horizon optimal feedback control for systems of the form \eqref{eq: dynamics} is linked to a state-dependent algebraic Riccati equation (ARE):
\begin{equation}\label{eq: sdre1}
A^\top(x)\Pi(x)+\Pi(x) A(x)-\Pi(x)B(x)R^{-1}B^\top(x)\Pi(x)+Q=0.
\end{equation}
We note that equation \eqref{eq: sdre1} is an ARE that changes every iteration, in fact it depends on the current state $x$. This makes the difference with respect to the standard LQR\footnote{LQR deals with the constant matrices $A(x(t))=A, B(x(t))=B$ in \eqref{eq: dynamics}.} problem where the ARE is constant and it is solved just once. SDRE might be thought as an MPC algorithm (see e.g. \cite{larsnmpc}) where the inner optimization problem is solved by \eqref{eq: sdre1}.

When equation \eqref{eq: sdre1} admits solution, it leads 
to a state-dependent Riccati operator $\Pi(x)$, from where we obtain a nonlinear feedback law given by
\begin{equation}\label{eq: sdref1}
u(x):=-R^{-1}B^\top(x)\Pi(x)x.
\end{equation}

We will refer to the feedback gain matrix as $K(x):=R^{-1}B^\top(x)\Pi(x).$
It is important to observe that the ARE \eqref{eq: sdre1} admits an analytical solution in a limited number of cases and the obtained control is only suboptimal. More general approaches following e.g. the dynamic programming approach \cite{B61} might be used. This goes beyond the scopes of this work, however one can easily replace, throughout the paper, the SDRE approach with a different (feedback) control method.

In \cite{BLT07}, it is shown that the SDRE method provides asymptotic stability if
	 $A(\cdot)$ is $\cC^1$ for $\|x\|\leq\delta$ and some $\delta >0$, $B(\cdot)$ is continuous and the pair $(A(x),B(x))$ is stabilizable for every $x$ in a non-empty neighbourhood of the origin. Thus, the closed-loop dynamics generated by the feedback law \eqref{eq: sdref1} are locally asymptotically stable. 
The SDRE algorithm proposed in \cite{BLT07} is summarized below.
\begin{algorithm}[ht]
	\caption{SDRE method }
	\label{alg: sdre}
 \begin{algorithmic}[1]
		\REQUIRE $\{t_0,t_1,\ldots\},$ model \eqref{eq: dynamics}, $R$ and $Q$,
		\FOR{$i=0,1,\ldots$}
		\STATE Compute $\Pi(x(t_i))$ from \eqref{eq: sdre1}
		\STATE Set $K(x(t_i)) := R^{-1}B^\top(x(t_i))\Pi(x(t_i))$
		\STATE Set $u(t):=-K(x(t_i))x(t), \quad \mbox{for }t \in [t_i, t_{i+1}]$ 
		\STATE Integrate the system dynamics with $u(t):=-K(x(t_i))x(t)$ to obtain $x(t_{i+1})$
		\ENDFOR
	\end{algorithmic}
\end{algorithm}


Assuming the stabilization hypothesis above, the main bottleneck in the implementation of Algorithm \ref{alg: sdre} is the high rate of calls to an ARE solver for \eqref{eq: sdre1}. We refer to \cite{AKS23} for efficient methods related to large scale problems. In this work, we will deal with small scale problems and thus solve the AREs with the Matlab function {\tt icare}.

\section{Identification and control of unknown nonlinear dynamics}\label{sec:idcon}

The system we want to identify and control is taken from \eqref{eq: dynamics} and reads
\begin{align}\label{eq:dynmu}
\begin{aligned}
	\dot x(t)&=\sum_{j=1}^n \mu_j A_j(x(t))x(t)+B(x(t))u(t),\quad t\in(0,\infty),\\ x(0)&=x_0,
	\end{aligned}
\end{align}
with the matrix function $A(x)$ in \eqref{eq: dynamics} given by  $ A(x) = \sum_{j=1}^n \mu_j A_j(x)$
and $A_j(x):\R^d\rightarrow\R^{d\times d}$ for $j= 1,\ldots,n$. The functions $A_j(x)$'s  may be thought as a library with terms that have to be selected by the coefficients $\mu_j$'s.
Note that this sum is not unique and that we
can include extra basis functions by simply setting the corresponding $\mu_j$'s to be zero. Throughout this work, we will assume that the terms $A_j(x)$'s and $B(x)$ are known.
Thus, the system \eqref{eq:dynmu} is fully identified by the knowledge of the coefficient $\mu=(\mu_1,\ldots,\mu_n)\in\R^n$ which is considered unknown in the present work.

We assume that there exists a true system configuration $
\mu^*\in\R^n$ which is not known but observable through the dynamics \eqref{eq:dynmu} setting $\mu=\mu^*$.
In other words, we assume that the dynamics generated by this true model configuration $\mu^*$ is always observable as a black box, i.e. if we choose a control we can compute the solution of \eqref{eq:dynmu} with the true parameter without knowing $\mu^*$ explicitly. This is a typical assumption in the Reinforcement Learning setting, where an agent can take actions and observe how the environment responds to them.\\

\noindent
 The cost functional we want to minimize is adapted from \eqref{disc_cost} and reads
\begin{equation}\label{disc_costmu}
\underset{u(\cdot)\in \mathcal{U}}{\min}\;J(u(\cdot),\mu^*)):=\int\limits_0^\infty \Big(\norm{x(t;\mu^*)}_Q^2+\norm{u(t)}_R^2\Big)\, dt,
\end{equation}
where the dependence on $\mu^*$ stresses that trajectories $x$ are observed from the true system configuration for a given input.


This addresses the problem of system-identification together with the control of \eqref{eq:dynmu}.
Indeed, we consider two unknowns: (i) the parameter configuration $\mu$ which is required to converge to $\mu^*$ and (ii) the control $u(t).$  The computation of the control will be done using Algorithm \ref{alg: sdre} which clearly depends on the parameter configuration. For an estimated parameter $\wt\mu\in\R^n$ such that $\wt\mu\neq \mu^*$,  the control will be denoted by $u(t;\wt \mu)$ to stress the dependence on the particular parameter configuration $\wt\mu$ in \eqref{eq:dynmu}. Instead the observed trajectory will be denoted by $x(t; u(t;\wt \mu),\mu^*)$. This notation considers a trajectory computed with the control $u(t;\wt \mu)$ plugged into the true system configuration. Furthermore, if we want to represent the solution at discrete time $t_i$ we will identify $x^i(u(t_i;\wt\mu),\mu^*) = x(t_i;u(t_i;\wt \mu),\mu^*).$

\begin{remark}[Notation]
Let us briefly summarize the notations valid throughout the whole paper for the parameter configuration: $\mu$ is a generic parameter, $\mu^*$ is the true system configuration and $\wt\mu$ is an estimate system configuration. 
\end{remark}

\subsection{The method}
\label{subsec:The method}
Let us now explain how we identify the system. We remind that our goal is to steer to the origin a partially unknown nonlinear system. We also aim at discovering the system on the fly through real-time observation of the trajectories. 
The workflow of our proposed method goes under the paragidm ``control first and identify later'' as follows:
\begin{enumerate}
    \item Pick a parameter configuration,
    \item Compute the corresponding  control,
    \item Observe the trajectories,
    \item Update the parameter configuration based on the observations, 
    \item Go to the second item.
\end{enumerate}


We use a Bayesian Linear Regression algorithm to estimate the system configuration $\mu^*$ from the observed trajectory data, as described in Section~\ref{sec:BLR}. We will now provide all the details of the proposed method which is summarized at the end of this subsection by Algorithm \ref{alg:RL}.

\paragraph{\bf Initial configuration}
To begin with we provide an initial estimate\footnote{The notation $\wt{\mu}^0$ refers to the parameter estimate at time $t_0$. This will become clearer later in this subsection} $\wt{\mu}^0\in\R^n$\ for the true system configuration $\mu^*$. To give an example, in the numerical tests we will set $(\wt{\mu}^0)_k\ = 1$ for $k=1,\ldots n,$ but if prior information about $\mu^*$ is available, it can be used to choose a proper $\wt{\mu}^0$. Note that $\wt\mu^0$ will act as $m_0$ in the prior distribution \eqref{dist:theta} of the BLR algorithm. We also need to choose an initial covariance matrix $\Sigma_0 \in \mathbb{R}^{n \times n}$. We observed heuristically that $\Sigma_0 = c I_n$, where $I_n$ is the $n\times n$ identity matrix and $c>0$ is large enough, works well in general.

\paragraph{\bf Computation of the control}
 At time $t_i$, we compute an approximate solution for the Algebraic Riccati Equation \eqref{eq: sdre1} corresponding to the current parameter estimate $ \wt{\mu}^i$. Then, we can set the feedback gain matrix $K(x(t_i); \wt{\mu}^i)$ and the feedback control $u(t;\wt{\mu}^i)$. 
 

\paragraph{\bf Observation of the trajectories}
At each iteration, we apply a piecewise-constant control $u(t;\wt{\mu}^i)$ with $t\in [t_i, t_{i+1}]$ and \textit{observe} the trajectory at time $t_{i+1}$, which will be either the actual trajectory, if we are dealing with a real physical system, or a simulated one, if we are simulating the physical system with some numerical methods. Thus, for a given configuration estimate $\wt\mu^i$ and its  control input $u(t_i;\wt\mu^i)$ computed following Algorithm \ref{alg: sdre}, we will observe the trajectory $x^{i+1}(u(t_i;\wt\mu^i);\mu^*)$.



The observation of the true trajectory has to be thought as a black box that provides the solution, or approximate solution, of the original controlled problem for a given control. The black box takes the control and the initial state as input and provides the trajectories as output. The observation of only one trajectory is due to the fact we aim at discovering and controlling on the fly updating the parameter estimate at each time instance. In RL, such methods are referred to as \textit{online}, in contrast with \textit{offline} methods where an agent can use multiple offline observations of the system to build a control.




\paragraph{\bf Update of the parameter estimate}
We now provide the crucial part of the method, that is how we update the parameter estimate using a Bayesian Linear Regression. To apply BLR and to obtain a problem in the form \eqref{eq: LR} we have to discretize the system \eqref{eq:dynmu}.
We provide, without loss of generality, an example through an implicit Euler scheme to discretize \eqref{eq:dynmu}.
Thus, the discretization of \eqref{eq:dynmu}, using e.g. an implicit Euler method and the correspondent feedback gain matrix $K^i:=K(x(t_i))$ (see Section \ref{sec: sdre}), reads
\begin{equation}\label{eq: discrete control matrix}
   \frac{x^{i+1}-x^i}{\D t} \approx  \sum_{j=1}^n \wt\mu_j^iA_j(x^{i+1})x^{i+1} - B K^i x^{i+1}, \quad i = 0,1,\ldots 
\end{equation}

where we have dropped the dependence on the control for $x$ and we recall that $x^i$ is the short notation for $x^i(u(t_i,\wt\mu^{i-1}),\mu^*)$. In our numerical simulations we deal with the control explicitly, that is why we use $K^i$ in \eqref{eq: discrete control matrix} and we consider the control piecewise constant in each interval $[t_i,t_{i+i}]$. We employ an implicit scheme due to numerical stability and our application to PDEs later in Section \ref{sec:test}.

Once we plug the true, observed trajectory in equation~\eqref{eq: discrete control matrix}, i.e. $x^i$ and $x^{i+1}$, we obtain a linear system of equations that the true system configuration $\mu_j^*$ solves, at least up to a certain approximation error. We use it to update our estimate $\wt\mu^i$ of the system configuration. Starting from equation~\eqref{eq: discrete control matrix}, we can write 
\[ \frac{x^{i+1}-x^i}{\D t} + B K^i x^{i+1} \approx \sum_{j=1}^n \wt\mu_j^i  A_j(x^{i+1})x^{i+1} . \]

\noindent
Thus, we obtain $d$ equations for the $n$ coefficients $\wt\mu_j^i$ as in \eqref{eq: LR}, which we can write in a more compact form
\begin{equation}\label{eq: alpha}
 Y^i \approx  X^i\wt\mu^i,
 \end{equation}
 where $\wt\mu^i\in\R^n, $ $X^i := [A_1(x^{i+1})x^{i+1},\ldots A_n(x^{i+1})x^{i+1}] \in \mathbb{R}^{d\times n}$ and $Y^i := \frac{x^{i+1}-x^i}{\D t} + B K^i x^{i+1} \in \mathbb{R}^d$.  
 The notation $\wt\mu^i$ stresses the fact that, at each time iteration we look for a parameter configuration that may differ on time.
 This problem fits into the structure presented in Section \ref{sec:BLR} and the solution for $\wt\mu^i$ is given by \eqref{thetaBLR}.



\paragraph{\bf Algorithm}Our proposed idea is finally summarized in Algorithm \ref{alg:RL} below.
\begin{algorithm}[H]
	\caption{Online identification and control}
	\label{alg:RL}
	\begin{algorithmic}[1]
		\REQUIRE $\{t_0,t_1,\ldots\},$ model $\{A_j(x)\}_{j=1}^n, B, R,Q$, $\wt{\mu}_0$,
  $\Sigma_0$
		\FOR{$i=0,1,\ldots$}
		\STATE Solve \eqref{eq: sdre1} and obtain $\Pi(x(t_i);\wt{\mu}^i)$ from \eqref{eq: sdre1} 
		\STATE Set $K(x(t_i);\wt{\mu}^i) := R^{-1}B^\top(x(t_i))\Pi(x(t_i);\wt{\mu}^i)$
		\STATE Set $u(t_i;\wt{\mu}^i):=-K(x(t_i);\wt{\mu}^i)x(t)$
		\STATE Apply the control $u(t_i;\wt{\mu}^i)$ and observe the trajectories $x^{i+1}(u(t_i;\wt \mu^i),\mu^*)$
        \STATE Compute $\wt{\mu}^{i+1}$ as in \eqref{thetaBLR} from \eqref{eq: alpha}
		\ENDFOR
	\end{algorithmic}
\end{algorithm}

\begin{remark}
    There might be cases in Algorithm \ref{alg:RL} where the ARE does not provide a solution and this will depend on the  parameter estimate. In that cases, we fix the feedback gain equal to the zero vector and we go to the next step.
\end{remark}

\begin{remark}\label{rem:3}
    Theoretical convergence of the parameter $\wt\mu^{i}$ to the true parameter $\mu^*$ for $i\rightarrow +\infty$ is not guaranteed and goes beyond the scopes of this paper. We decided to keep this study for a follow up work. However, in the numerical tests in Section \ref{sec:test}, we observed numerical convergence of the method. The identification of the system configuration can be stopped if, for a certain $\bar i>0$, we obtain $\|\wt{\mu}^{\bar{i}}-\wt{\mu}^{\bar{i}-1}\|_\infty<tol_\mu$ with $tol_\mu>0$ being the desired threshold.  Note this criteria is only heuristic and that Algorithm \ref{alg:RL}, as it is, does not need the parameter convergence. The primary goal is to stabilize an unknown  control system at $0$.
    \end{remark}

\subsection{Application to PDEs}
Our ultimate goal is the application of Algorithm \ref{alg:RL} to identify and control nonlinear PDEs given by
\begin{equation}\label{eq:pde}
    \left\{ \begin{aligned}
    y_t(t,\xi) &= \sum_{j=1}^n \mu_j F_j(y(t,\xi),y_\xi(t,\xi),y_{\xi\xi}(t,\xi),y_{\xi\xi\xi}(t,\xi),\ldots) + B^T(\xi) u(t), \qquad &&t \in [0,+\infty), \xi \in (a,b), \\
    y(0,\xi) &= y_0(\xi), \qquad &&\xi \in [a,b], \\
    y(t,a) &= 0, \ y(t,b) = 0, \qquad &&t \in [0,+\infty).
    \end{aligned} \right.
\end{equation}
where $y:[0,\infty]\times\R\rightarrow\R$, $\mu_j\in\R,$ $u(t):[0,\infty)\rightarrow\R^m$ and 
$B(\xi):[a,b]\rightarrow\in\R^m$. Without loss of generality, we set zero Dirichlet boundary conditions.
 We assume that the model is given by the sum of simple monomial bases functions $F_j$ of $y$ and its derivatives. Similarly to \eqref{eq:dynmu}, the functions $F_j$'s  may be thought as a library with terms that has to be selected by the coefficients $\mu_j$'s. Note that this sum is not unique and that we can include extra basis functions by simply setting the corresponding $\mu_j$'s to be zero. 

The numerical discretization of \eqref{eq:pde}, by e.g. finite differences method \cite{L07}, provides a system in the form \eqref{eq:dynmu}, where each component of $x\in\R^d$ corresponds to the grid points, say $x_i(t)\approx y(t,\xi_i)$ for $i=1,\ldots, d$ and $A_j(x(t))x(t)\approx F_j$, where $F_j\in \mathbb{R}^d$ denotes the 
the basis function evaluated at all the grid points such that $(F_j)_i=F_j(y(t,\xi_i),y_{\xi}(t,\xi_i),y_{\xi\xi}(t,\xi_i),y_{\xi\xi\xi}(t,\xi_i),\ldots).$  In Section \ref{sec:test}, we will explain in detail how to obtain each term $A_j(x)$. We note that the matrices $A_j(x)$ take into account the boundary conditions.




The continuous cost functional we want to minimize is 
\begin{equation}\label{cost_func}
    J(u;\mu^*) = \int_0^\infty \left(\|y(t,\cdot;\mu^*)\|_{L^2(a,b)}^2 + \|u(t)\|_R^2\right)\,dt
\end{equation}
with $R$ defined after equation \eqref{eq: dynamics} and we stress the dependence of the trajectory $y$ on the true system configuration. The discretization of \eqref{cost_func} corresponds to the choice $Q=\Delta \xi I_d$ in \eqref{disc_cost} with $\Delta \xi>0$ being the spatial step size and $I_d$ the $d\times d$ identity matrix.



\section{Numerical experiments}\label{sec:test}

In this section, we will show our numerical examples to validate the proposed method. To set the section into perspective we provide the continuous PDE model studied which reads:
\begin{equation}\label{eq: parametric PDE}
    \left\{ \begin{aligned}
    y_t(t,\xi) &= \mu_1 y_{\xi\xi}(t,\xi) + \mu_2 y_\xi(t,\xi) +\mu_3 y(t,\xi) + \mu_4 y^2(t,\xi) \\
    &   + \mu_5 y^3(t,\xi) + \mu_6 y (t,\xi) y_\xi(t,\xi) +\mu_7 y_{\xi\xi\xi}(t,\xi) + B^T u(t) \quad &&t \in [0,t_{end}], \xi \in (a,b) \\
    y(0,\xi) &= y_0(\xi) \quad &&\xi \in [a,b] \\
    y(t,a) &= 0, \ y(t,b) = 0  \quad &&t \in [0,t_{end}] .
    \end{aligned} \right.
\end{equation}
where for numerical reasons we have to choose a finite horizon with $t_{end}>0$ large enough to simulate the infinite horizon problem and such that for $t>t_{end}$ the controlled solution will not change significantly. We remark that, based on the choice of the parameters, the model \eqref{eq: parametric PDE} includes the control of e.g. heat equation, advection equation, diffusion-reaction-convection equation, burgers equation, viscous burgers, etc. Many of these models have different physical interpretations between them.

In this model we fix $n=7$ libraries and
in order to fit into the desired canonical form \eqref{eq:dynmu} we use Finite difference (FD) method (see e.g. \cite{L07}) where the discrete state $x(t)$ corresponds to the approximation of $y(t,\xi)$ at the grid points. 

The term $A(x)$ will be given by
$$A(x) = \mu_1 \Delta_d + \mu_2 T + \mu_3 I_d + \mu_4 diag(x) +\mu_5 diag(x\circ x)+\mu_6 \tilde{D}(x) +\mu_7 M$$
where the symbol $\circ$ denotes the Hadamard or component-wise product and

\begin{itemize}
\item $\Delta_d\in\R^{d\times d}$ is the FD approximation of the Dirichlet Laplacian with $\Delta_d:=\Delta \xi^{-2}{\tt tridiag}([1\ -2 \ 1],d),$\footnote{The notation {\tt tridiag}([a b c],d) stands for a tridiagonal $d\times d$ matrix having the constant values $b\in\R$ on the main diagonal, $a\in\R$ on the lower diagonal and $c\in\R$ on the upper diagonal.}
\item $T\in\R^{d\times d}$ is the FD upwind or downwind approximation of the advection term such that
\begin{itemize}
    \item if $\mu_2>0$, $T = T_{neg} = \Delta \xi^{-1}{\tt tridiag}([-1\ 1\ 0],d)$, 
    \item if $\mu_2<0$, $T = T_{pos} = \Delta \xi^{-1}{\tt tridiag}([0\ -1\ 1],d),$
\end{itemize} 
\item $I_d\in\R^{d\times d}$ is the identity matrix,
\item $diag(x)\in\R^{d\times d}$ is a diagonal matrix with the components of the vector $x$
\item $\tilde D(x)\in\R^{d\times d}$ is a matrix such that its $i$-th row $\tilde D(x)_i$ is 
\begin{itemize}
    \item if $(\mu_6x)_i>0$, $\tilde D(x)_i=(diag(x) T_{neg})_i$
    \item if $(\mu_6x)_i<0$, $\tilde D(x)_i=(diag(x) T_{pos})_i$
\end{itemize}
where $(\mu_6x)_i$ indicates the $i$-th element of the vector $\mu_6x$ and $(\tilde D(x))_i$ indicates the $i$-th row of the matrix between parentheses
\item $M\in\R^{d\times d}$ is a FD approximation of the third order derivative: $M=-\frac{1}{2\Delta \xi^3}{\tt pentadiag}([1 \ -2 \ \ 0 \ \ 2 \ -1],d)$
\footnote{The notation {\tt pentadiag}([a b c e f],d) stands for a pentadiagonal $d\times d$ matrix having the constant values $c\in\R$ on the main diagonal, $b\in\R$ on the lower and $e\in\R$ on the upper diagonal and $a\in\R$ on the second diagonal below and $f\in\R$ on the second diagonal above the main diagonal.}.
\end{itemize}
Controlled trajectories are integrated in time using an implicit Euler method (see \eqref{eq: discrete control matrix}), which is accelerated using a Jacobian--Free Newton Krylov method (see e.g. \cite{KK04}) using $10^{-5}$ as threshold for the stopping criteria of the Netwon method and less than $500$ iterations.
As mentioned in Remark \ref{rem:3}, in our numerical simulations we have observed convergence of our estimated configuration to the true one. Therefore, we have added to Algorithm \ref{alg:RL} the stopping criteria with $tol_\mu=10^{-5}$.


We present three numerical test cases with nonlinear PDEs. Those are the PDEs we can observe.
The first test is a nonlinear diffusion-reaction equation, known as the Allen-Cahn equation ($\mu^*=$[1, 0, $11$, 0, $-11$, 0, 0]). The second test studies the viscous Burgers' equation ($\mu^*=$[0.01, 0, 0, 0, 0, 1, 0]) and the third one the so called Korteweg-De Vries (KdV) model ($\mu^*=$[$0.5$, 0, 0, 0, 0, 6, $-1$]). The goal of all our tests is the stabilization of the (unknown) dynamics to the origin by means of the minimization of the cost functional \eqref{cost_func}, that can be approximated $\|y(t,\cdot;\mu)\|_{L^2(a,b)}^2 \simeq \sum_{i=1}^d \Delta \xi y(t,\xi_i;\mu)^2 =x^T(t)Qx(t)=\norm {x(t)}^2_Q$, where $\xi_i = a + \Delta \xi i$ and $Q=\Delta \xi I_d$ and with $R=0.01$, thus obtaining \eqref{disc_costmu}. 

In all our tests, we will plot (i) the \textit{uncontrolled solution} of the true dynamical system, i.e. the solution of \eqref{eq: parametric PDE} obtained choosing $u(t)\equiv 0$ and $\mu = \mu^*$, (ii) the \textit{controlled solution} based on the SDRE method where $\mu^*$ is known, and (iii) our \textit{RL solution} identified by Algorithm \ref{alg:RL} where $\mu^*$ has to be discovered. We will then compare the optimal control computed by Algorithm \ref{alg: sdre} and our method Algorithm \ref{alg:RL} and the evaluation of the cost functionals. Furthermore, the history of the estimated coefficients $\wt\mu^i$ over time will be presented.
Finally, we will also discuss a numerical convergence towards the control of the continuous PDE problem. To do that, we will compute the control for a given spatial discretization $\Delta \xi$ and show that the obtained control is robust enough to stabilize the same problem with decreasing values of the spatial discretization. This will show the numerical mesh independence and the robustness of the proposed method.

The RL assumption relies on the observability of the dynamics with the true system configuration $\mu^*$. This should be thought as a black box where the true model can be computed (or approximated). In this work, since we do not know the exact solution, we will use two different numerical approaches to obtain the observed trajectories: (i) we use the same scheme, e.g. backward Euler method, used in \eqref{eq: discrete control matrix} but with the true parameter $\mu^*$ and (ii) an explicit Runge Kutta scheme for stiff problems. 

For the sake of completeness, we provide some more numerical details on the two schemes. Again, we stress that those details are not critical to the algorithm, but they are only needed for the numerical simulations. Indeed, one could use any method or even a "real" black box.\footnote{By the term "real" black box we intend something that takes an input and provides the trajectories without knowing how they are computed.}
These two methods will have a different way to approach the feedback control. Indeed, the first approach is "implicit" in the control term $Kx$, so will be called "implicit approach" or "implicit algorithm" in the following; in this case, we will have the feedback control in the form  $K^{i}x^{i+1}$, mainly for stability reasons. The second scheme is explicit and the feedback control will be $K^i x^{i}.$ In the paper, the latter approach has been implemented using the Matlab function {\tt ode15s}. We remark that the second approach could be used in a real application, replacing the result of the {\tt ode15s} function with an observation of the system evolution, and will be called "black box algorithm". Note that in this case we have e.g. $Y^i=\frac{x^{i+1}-x^i}{\D t} + B K^i x^{i}$ in \eqref{eq: alpha}.

We remark that in both cases, say the use of the implicit algorithm or the use of the "black box", we added noise to the data used for regression. After computing the control $u(t_i)$ and the trajectory $x(t_{i+1})$, we obtained the matrix $X=[A_1(x(t_{i+1}))x(t_{i+1}), \ldots, A_n(x(t_{i+1}))x(t_{i+1})]\in \mathbb{R}^{d\times n}$. To each column of $X$ we added a vector of independent Gaussian random variables, each with mean 0 and standard deviation given by $0.01$ times the mean of the absolute values of the components in the column itself as follows: 
\begin{equation*}
    {X}_{i,j} \leftarrow {X}_{i,j} + \mathcal{N}\biggl(0,\Bigl(\frac{0.01}{d}\sum_{k=1}^d |X_{k,j}|\Bigr)^2\biggr).
\end{equation*}
This will be referred as $1\%$ relative noise in the following, and has been used in every numerical tests, in order to simulate noise on data from real applications. The noise can be also interpreted as a variation on the observed system that adds negligible terms not in the library.

Finally, we note that in all the tests the \emph{prior distribution} on the parameter $\mu$ was initialized as described in Section \ref{subsec:The method}, i.e. we started with a normal distribution with mean $\wt \mu^0=[1, 1, 1,  1,  1,  1,  1]^T$ 
and covariance matrix $\Sigma_0 = cI_7$ 
with $c=1000$ for test 2 and 3 and $c=200000$ for test 1. In general, $c>0$ must be chosen large enough to guarantee flexibility to the model. Indeed, the smaller it is, the closer the final approximation of $\mu$ will tend to be to the chosen initial \textit{prior} distribution.

The tests presented in this paper were performed on a DELL Latitude 7200, Intel(R) Core(TM) i5-8265U CPU 1.60GHz, using \textsc{Matlab}.

\subsection{Test 1: Allen-Cahn}
\begin{figure}[htbp]
\centering
\includegraphics[width=.3\textwidth]{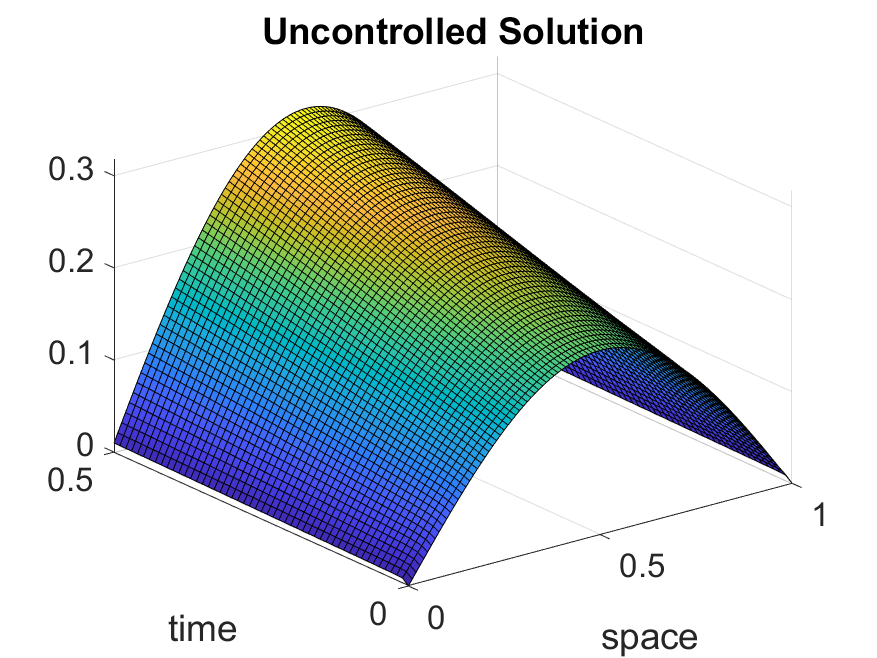}
\includegraphics[width=.3\textwidth]{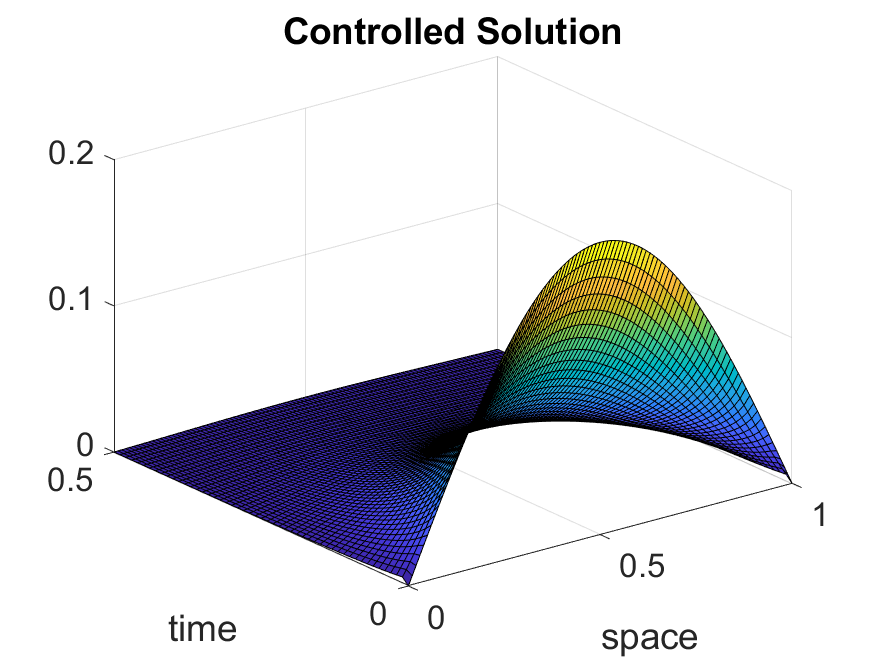}
\includegraphics[width=.3\textwidth]{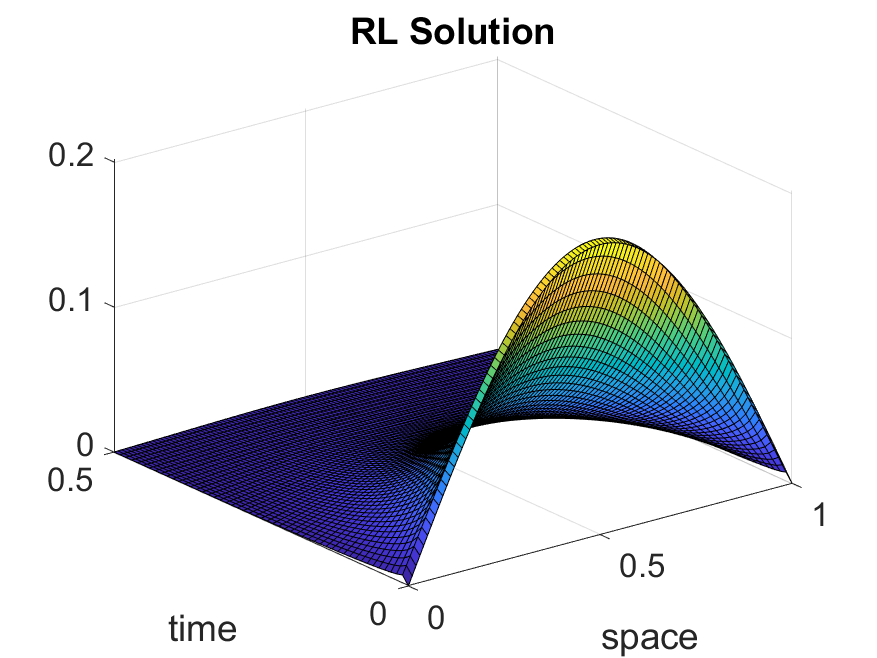}
\caption{Test 1: Allen-Cahn, $\Delta\xi = 0.01$, $\Delta t = 0.01$, 1\% relative noise.}
\label{fig: Allen-Cahn Grune}
\end{figure}

Our first test is inspired by the example in \cite{AGW10} where it is shown that a MPC approach, for short prediction horizion, does not stabilize the following equation:
\begin{equation*}
    \left\{ \begin{aligned}
    y_t(t,\xi) &= y_{\xi\xi}(t,\xi) + 11(y(t,\xi) - y^3(t,\xi)) + u(t), \qquad &&t \in (0,0.5], \ \xi \in (0,1)\\
    y(0,\xi) &=0.2\sin (\pi \xi) \qquad && \xi \in (0,1), \\
    y(t,0) &= 0, \ y(t,1) = 0, \qquad &&t \in [0,0.5].
    \end{aligned} \right.
\end{equation*}
This model is known as Allen-Cahn equation, or Chaffee-Infante equation. Note that for this example a small horizon with $t_{end}=0.5$ is enough to simulate the infinite horizon problem since the control problem will be stabilized before as it is shown in Figure \ref{fig: Allen-Cahn Grune}. Here, we use the same settings of \cite{AGW10}, i.e. $\Delta \xi=0.01=\Delta t$ and we obtain its discrete version as described in \eqref{eq:dynmu} where the $B$ vector is given by a vector of ones. The dimension of the discrete problem is $d = 101.$ The only difference with respect to \cite{AGW10} is that we 
introduce the control as time-dependent function instead of dealing with a control as a function of time and space. The parameter used in the observable trajectories are $\mu_1^* = 1, \mu_3^* = 11,\mu_5^*=-11$ in \eqref{eq: parametric PDE} subject to the cost functional recalled in Section \ref{sec:test}.
In the left panel of Figure~\ref{fig: Allen-Cahn Grune}, we show the solution of the uncontrolled problem whereas in the middle panel the trajectory computed using Algorithm \ref{alg: sdre}. Both simulations have been computed knowing the true system configuration. It is clear that the solution is stabilized.
We remark that the SDRE method is able to stabilize the problem with a an infinite prediction horizon for the linearized problem, whereas the method in \cite{AGW10} uses a finite prediction horizon, that is required to be of size at least $11\Delta t$, for the nonlinear equation. 
Clearly, our inner minimization problem is different from the approach proposed in \cite{AGW10} but nevertheless it is interesting to see its stabilization through SDRE.

Finally, in the right panel of Figure \ref{fig: Allen-Cahn Grune}, we show the solution of Algorithm \ref{alg:RL}. It is clear that  with our method we can also stabilize the problem and identify the correct model as shown in Table \ref{tab: Allen-Cahn Grune}.

In Table \ref{tab: Allen-Cahn Grune} we show the results of Algorithm \ref{alg:RL} concerning the paramater configuration $\wt\mu$ estimated. We can see that the reconstructed values (second row of the table) are very close to the desired configuration considering the discretization $\Delta \xi = 0.01 = \Delta t$ and the noise added at each iteration.
\begin{table}[htbp]
    \centering
        \begin{tabular}{r|r|r|r|r|r|r|r}
        True $\mu^*$ & 1 & 0 & 11 & 0 & -11 & 0 & 0\\ 
        \hline
        Estimated $\wt\mu$ & 0.9992 & -0.0017 & 11.0008 & -0.0653 & -10.8232 & 0.0431 & 0
    \end{tabular}
    \caption{Test 1: Reconstruced parameter configuration for Allen-Cahn with $\Delta \xi = 0.01$, $\Delta t = 0.01$, 1\% relative noise.} 
\label{tab: Allen-Cahn Grune}
\end{table}

\begin{figure}[htp]
\centering
\includegraphics[width=.3\textwidth]{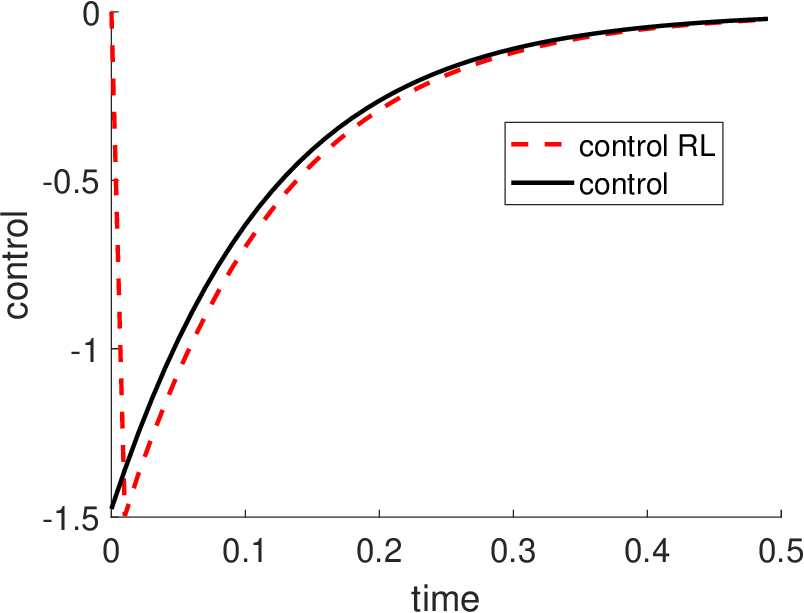}\hfill
\includegraphics[width=.3\textwidth]{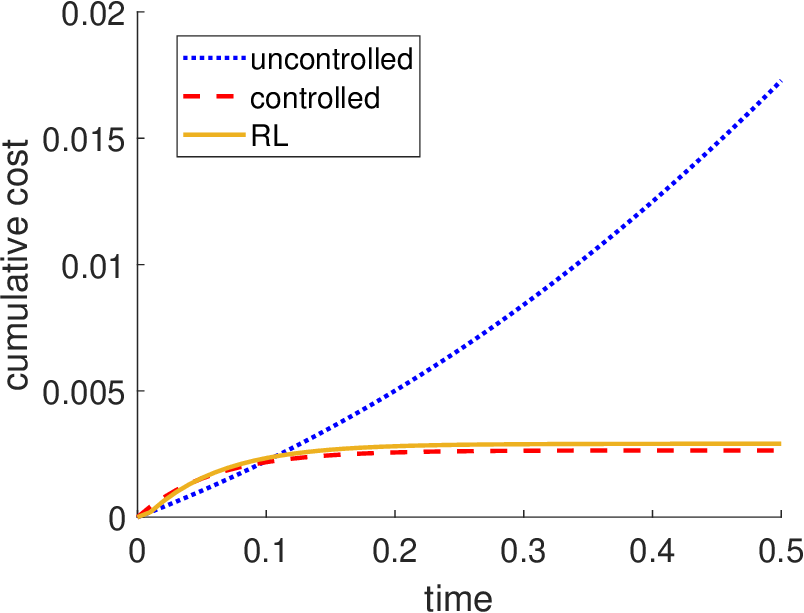}\hfill
\includegraphics[width=.3\textwidth]{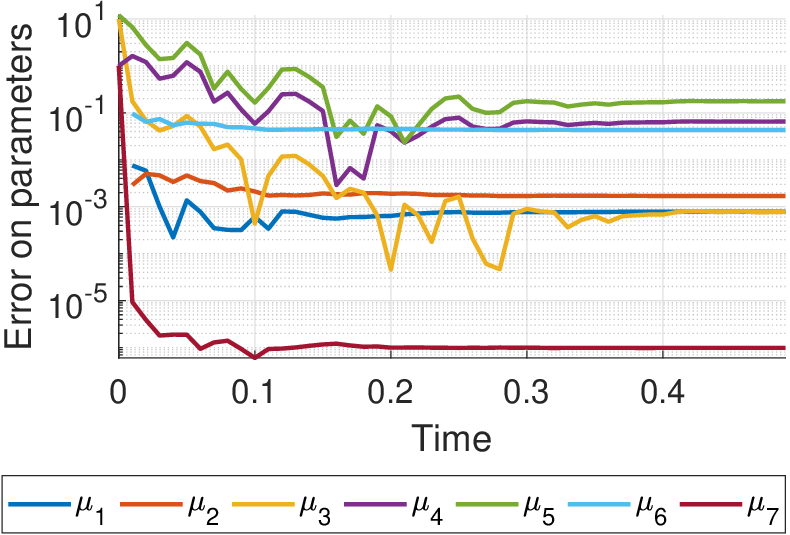}
\caption{Test 1: Allen-Cahn, $\Delta \xi = 0.01$, $\Delta t = 0.01$, 1\% relative noise. On the left, the comparison between the control found using knowledge of the true $\mu$ and the control found by the RL algorithm is shown. In the middle, the cumulative cost. On the right, the error on the parameter estimation at each time.}
\label{fig: Allen-Cahn Grune u mu}
\end{figure}

In Figure \ref{fig: Allen-Cahn Grune u mu}, we compare the control of the SDRE algorithm and the RL based one. One can see that at the beginning the RL control starts from 0 because we decided not to act at the first iteration since we do not have information at that stage. Then, we can see that, slowly, the RL control tends to the SDRE one which is our reference control. In the middle plot of Figure \ref{fig: Allen-Cahn Grune u mu}, we show the evaluation of the cost functional. One can see that the RL algorithm is very close to  the SDRE method and, as expected, the SDRE cost functional provides lower values. This is clear since our method starts without any knowledge of the model which is learnt on the fly. Finally, for completeness, in the right panel of Figure \ref{fig: Allen-Cahn Grune u mu} we show the convergence history of the parameter configurations. In this example, until the end of the chosen time interval, the algorithm never stops updating the distribution. It would stop at t=0.58 (after 58 iterations) if a longer time interval was considered.

\subsection{Test 2: Viscous Burgers }
The equation we study in this the test is the viscous Burgers problem which reads:
\begin{equation}
\label{eq:burg}
\left\{ \begin{aligned}
    y_t(t,\xi) &= 0.01 y_{\xi\xi}(t,\xi) + y(t,\xi) y_\xi(t,\xi) + 
   B(\xi)^T u(t),  &&t \in [0,2], \ \xi \in (-1.5,1.5),\\
    y(0,\xi) &=\sin (\pi \xi)\chi_{[0,1]}(\xi)  &&\xi \in (-1.5,1.5),\\
    y(t,-1.5) &= 0= y(t,1.5),  &&t \in [0,2].
    \end{aligned} \right.
\end{equation}
with $B(\xi)^T =  \begin{pmatrix}
        \chi_{[0.25,0.5]}(\xi),  & \chi_{[0.75,1]}(\xi)
    \end{pmatrix} $. The true system configuration is given by $\mu_1^* = 0.01$ and $\mu_6^*=1$ in \eqref{eq: parametric PDE}. We note that in this example $u(t)\in\R^2$, and we set in \eqref{disc_costmu} $Q=\Delta \xi I_d, R =0.01$. The discretization of \eqref{eq:burg} is done with $\Delta \xi = 0.025 = \Delta t.$ This discretization leads to a problem of dimension $d=121$.
    In Table~\ref{tab:Viscous_Burger}, one can find the true coefficients versus the reconstructed ones at the last iteration using Algorithm \ref{alg:RL}. We can see that our algorithm matches the desired configuration considering the order of the finite discretization used.


\begin{table}[thb]

\centering
\begin{tabular}{c|r|r|r|r|r|r|r}
True $\mu^*$ & 0.01 & 0 & 0 & 0 & 0 & 1 & 0\\ 
\hline
Estimated $\wt\mu$ & 0.0096 & 0 & -0.0008 & 0.002 & -0.001 & 0.9999 & 0
\end{tabular}
\caption{Test 2. Reconstructed parameter configuration for viscous Burger with $\Delta \xi = 0.025$, $\Delta t = 0.025$, 1\% relative noise.} 
\label{tab:Viscous_Burger}
\end{table}



\begin{figure}[htbp]
\centering
\includegraphics[width=.32\textwidth]{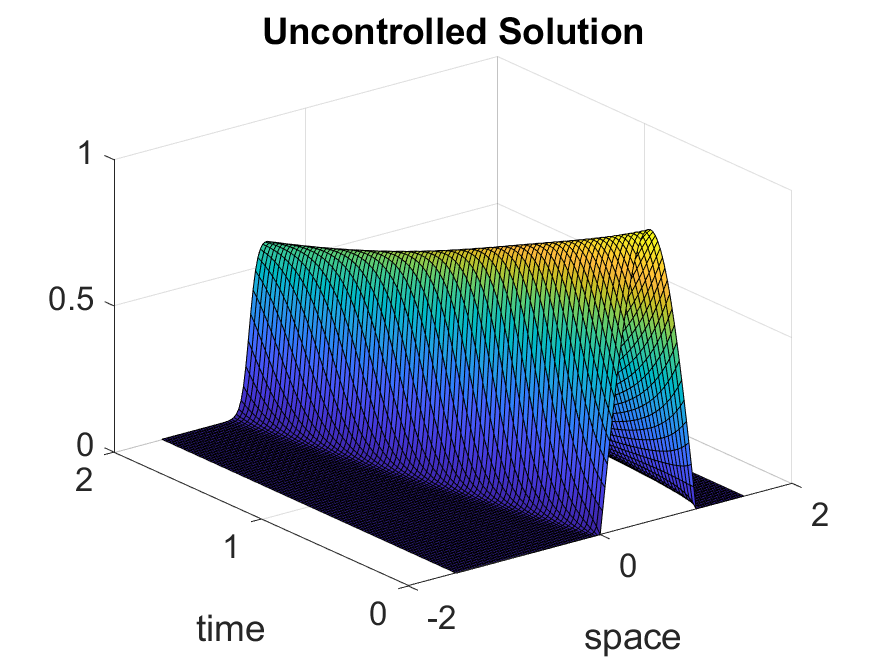}
\includegraphics[width=.32\textwidth]{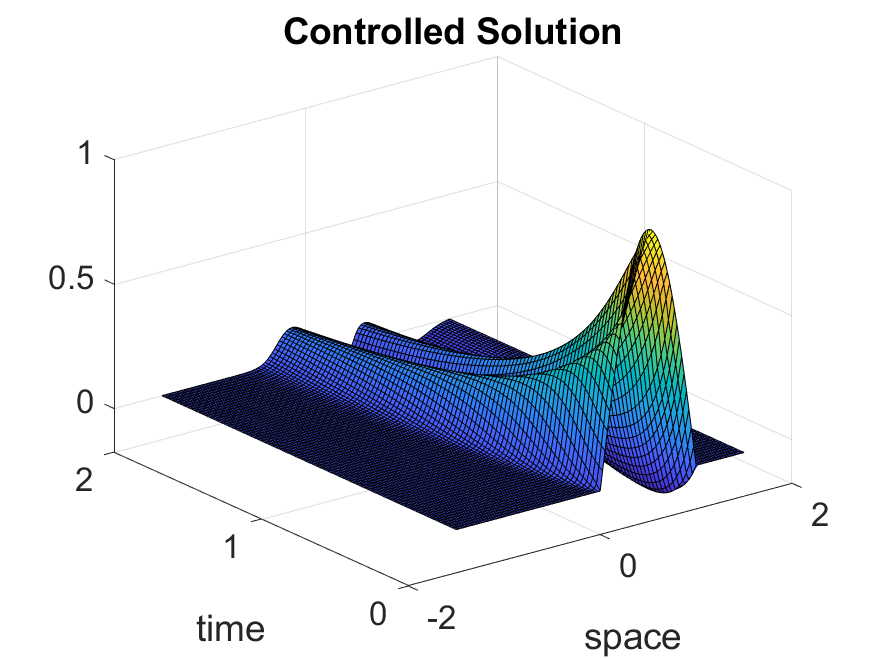}
\includegraphics[width=.32\textwidth]{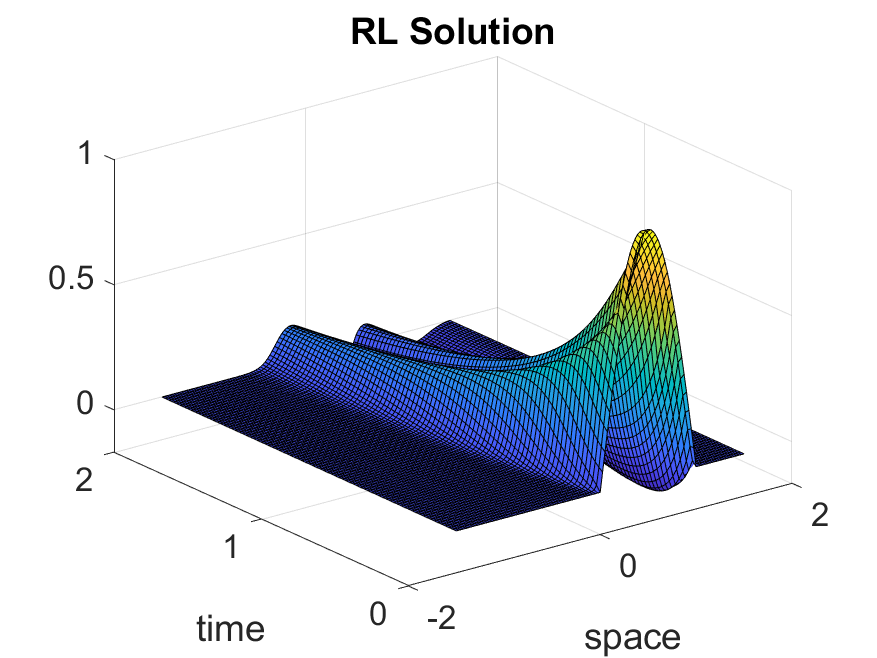}
\caption{Test 2: Viscous Burgers, $\Delta \xi = 0.025$, $\Delta t = 0.025$, 1\% relative noise.}
\label{fig: Viscous Burgers}
\end{figure}

The three trajectories are compared in Figure \ref{fig: Viscous Burgers}. One can see a clear difference between the uncontrolled solution (left panel) and the other two plots. Moreover the controlled trajectories show a similar behavior.

\begin{figure}[htp]
\centering
\includegraphics[width=.3\textwidth]{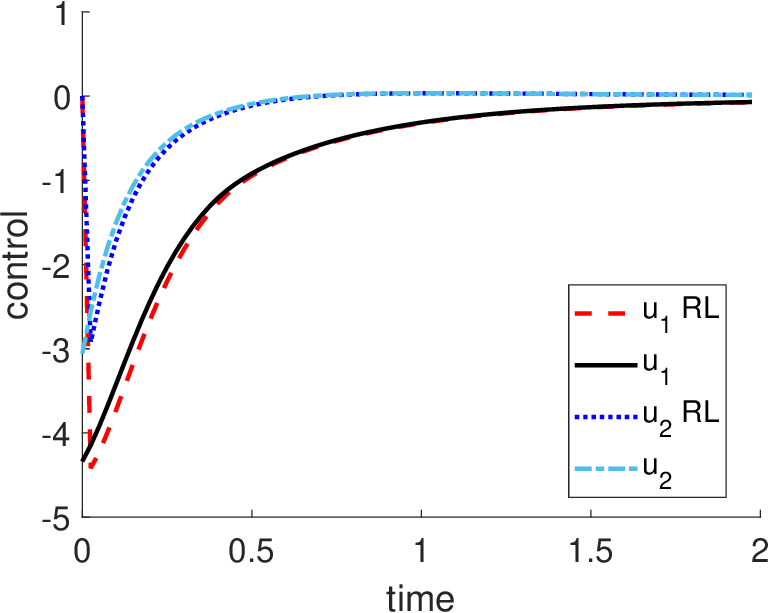}\hfill
\includegraphics[width=.3\textwidth]{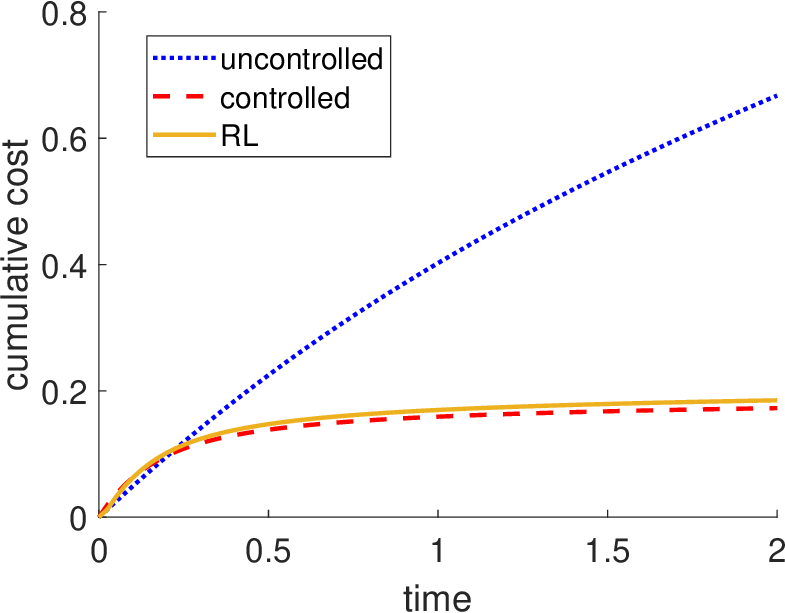}\hfill
\includegraphics[width=.3\textwidth]{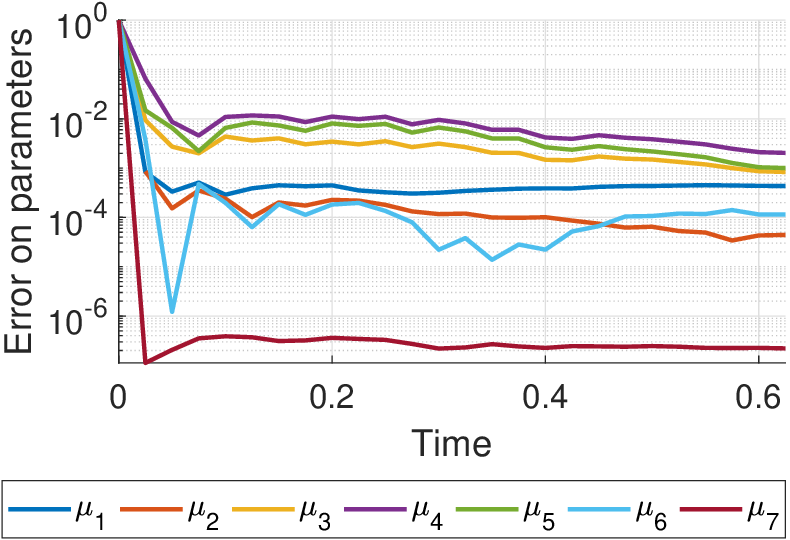}
\caption{Test 2: Viscous Burgers, $\Delta \xi = 0.025$, $\Delta t = 0.025$, 1\% relative noise. The first plot shows the comparison between each of the two components of the control found using knowledge of the true $\mu$ and the control found by the RL algorithm. The other plots show the cumulative cost and the error on the parameter estimation at each time  until the update stop.}
\label{fig: Viscous Burgers u mu}
\end{figure}

A more detailed comparison between Algorithm \ref{alg: sdre} and Algorithm \ref{alg:RL} is shown in Figure \ref{fig: Viscous Burgers u mu}. Indeed, in the left plot we show the two controls obtained from each algorithm that are very close to each other. The evaluation of cost functional is shown in the middle plot of Figure \ref{fig: Viscous Burgers u mu} and one can see that, again as expected, the value of the RL based method is slightly larger to the SDRE method. Finally, in the right plot we show the error in the convergence of the parameter configuration. The method stops updating the configuration estimate at time t=0.625 (i.e. after 25 iterations out of 80).

\begin{figure}[htbp]
\centering
\includegraphics[width=.4\textwidth]{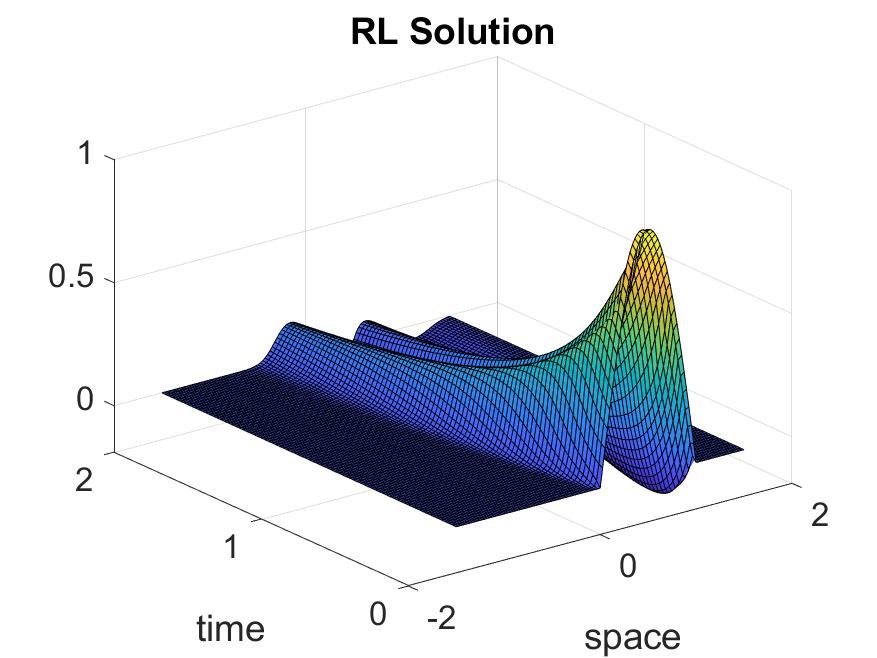}
\includegraphics[width=.4\textwidth]{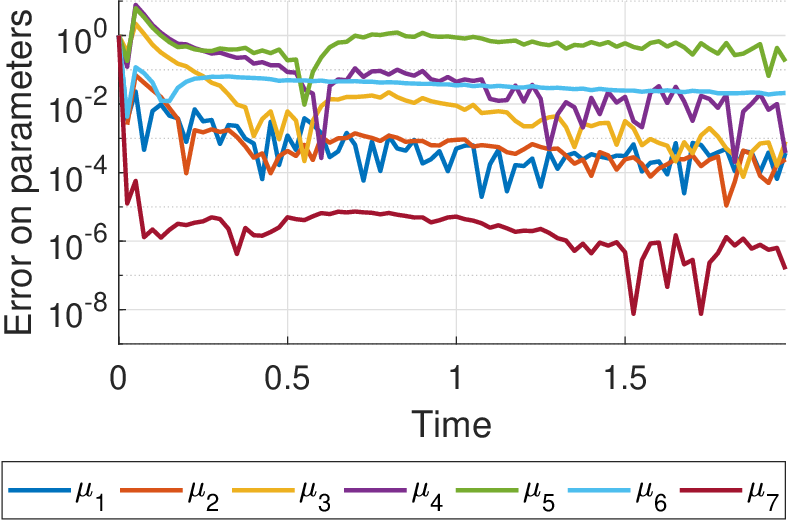}\\
\includegraphics[width=.4\textwidth]{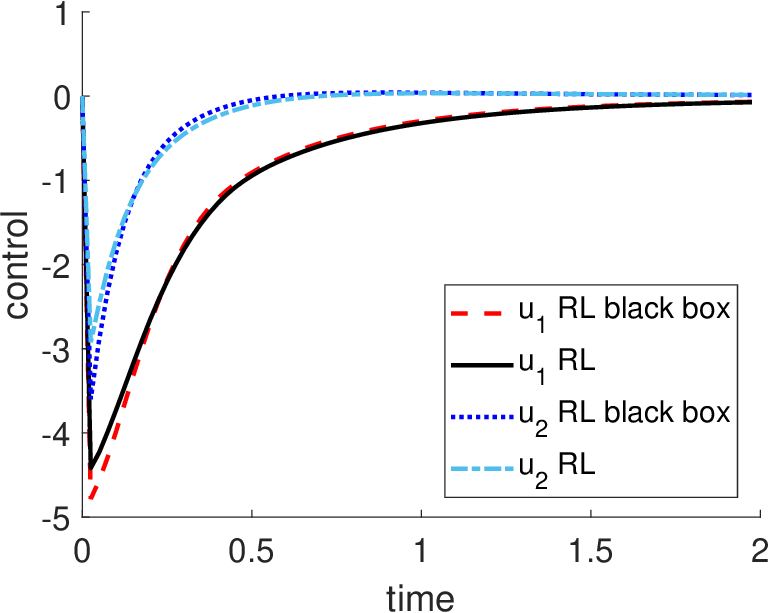}
\includegraphics[width=.4\textwidth]{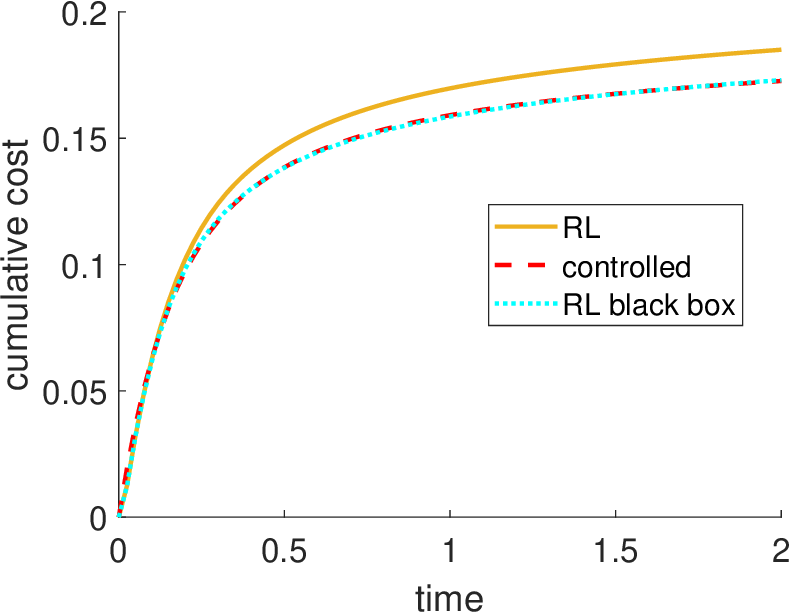}
\caption{Test 2: Viscous Burgers, $\Delta \xi = 0.025$, $\Delta t = 0.025$, 1\% relative noise. Results with a black box, all parameters considered.}
\label{fig: Viscous Burgers bb_full}
\end{figure}

\paragraph{\bf Results with a black box.}
For this test case we also show the results obtained using a real black box. At each iteration $t_i$ we first solved the Riccati equation, thus obtaining $K^i$, then we found the control $u=K^ix^i$ and finally we let the system evolve, thus obtaining $x^{i+1}$. For the evolution, we used the \textsc{Matlab} function {\tt ode15s} at each iteration. 
Table \ref{tab: Viscous Burgers bb_full} shows the final approximations of the parameters. We can see that the term $\wt\mu_5$, which was close to the correct value $0$ using an implicit scheme in Table \ref{tab:Viscous_Burger}, appears in the reconstruction with a value of ~0.18 in this case. Nevertheless, the algorithm is still able to provide a stabilizing control which is very close to the one obtained using the implicit scheme, as shown in Figure~\ref{fig: Viscous Burgers bb_full}.

\begin{table}[htbp]
\centering
\begin{tabular}{r|r|r|r|r|r|r|r}
True $\mu^*$ & 0.01 & 0 & 0 & 0 & 0 & 1 & 0\\ 
\hline
Estimated $\wt\mu$ & 0.0096 & 
-0.0002 & 0.0008 & 0.0004 & 0.1762 & 1.021 & 0
\end{tabular}
\caption{Test 2: Viscous Burgers, $\Delta \xi = 0.025$, $\Delta t = 0.025$, 1\% relative noise. Results with a black box, all parameters considered.} 
\label{tab: Viscous Burgers bb_full}
\end{table}

We then tested our algorithm using different libraries, i.e. considering only some $F_j$'s. In Table \ref{tab: Viscous Burgers bb_partial} we report the results when not considering the 5th term $y^3(t, x)$, whose parameter $\mu_5$ is the extra term appearing when working with the whole library. It is then clear that our method works accurately. Furthermore, we report in Table \ref{tab: Viscous Burgers bb_partial_16} the results when considering only the terms that belong to the problem; again, our algorithm was able to approximate them well.

\begin{table}[htbp]
\centering
\begin{tabular}{r|r|r|r|r|r|r|r}
True $\mu^*$ & 0.01 & 0 & 0 & 0 & $-$ & 1 & 0\\
\hline
Estimated $\wt\mu$ & 0.0101 & -0.0003 & 0.0007 & 0.0116 & $-$ & 1.0199 & 0 \end{tabular}
\caption{Test 2: $\Delta \xi = 0.025$, $\Delta t = 0.025$, 1\% relative noise. Results with a black box, $\mu_5$ not considered.} 
\label{tab: Viscous Burgers bb_partial}
\end{table}

\begin{table}[hthb]
\centering
\begin{tabular}{r|r|r|r|r|r|r|r}
True $\mu^*$ & 0.01 & - & - & - & - & - & 1\\
\hline
Estimated $\wt\mu$ & 0.0099 & - & - & - & - & - & 1.0564
\end{tabular}
\caption{Test 2, $\Delta \xi = 0.025$, $\Delta t = 0.025$, 1\% relative noise. Results with a black box, only $\mu_1$ and $\mu_6$ considered.} 
\label{tab: Viscous Burgers bb_partial_16}
\end{table}

The reason why the algorithm is not able to approximate $\mu_5$ might be that the components of the solution vector $x$ all tend to zero with the applied control, so the components of $x^3$ (each element of $x$ raised to the power of 3), which is the term that must be multiplied by $\mu_5$, tend to 0 very rapidly. This is also justified by the error in infinity norm we computed in Table \ref{tab: cfr bb}. There, we have computed the difference between the controlled solution and our RL solutions using the full library in the first column, the library without the $\mu_5-$term in the second column, and the library with only the correct terms in the last column. One can see that there is no difference in using the libraries chosen. Indeed, the computed controls appear to be the same and our method is always able to stabilize the problem even in the case of the full library (see Table \ref{tab: Viscous Burgers bb_full}). This further validates our method which is able to stabilize the problem even if the discovered model does not match perfectly with the true system configuration. The reason is that even if a configuration $\wt\mu$ doesn't match exactly the true configuration $\mu^*$, by construction it solves the linear system \eqref{eq: discrete control matrix} and so it well approximates the system dynamics, at least along the controlled trajectory.



\begin{table}[t]
\centering
\begin{tabular}{c|c|c|c}
 & full library & No $\mu_5$ & Only $\mu_1,\mu_6$\\
\hline
$\|y_{RL}- y_{RL-bb}\|_{\infty}$ & 0.054254 & 0.054254 & 0.054254\\
 \hline
$\|y_c-y_{RL-bb}\|_{\infty}$ & 0.113177 & 0.113177 & 0.113177\\
\hline
$\|y_c-y_{RL}\|_{\infty}$ & 0.111801 & 0.111801 & 0.111801\\
\hline
\end{tabular}
\caption{Difference between the controlled approximation $y_c$ with Algorithm \ref{alg: sdre} for the known problem, the RL approximation $y_{RL}$ and the RL with a black box $y_{RL-bb}$. In the first column the error is computed using the results with a full library, in the second excluding the term $\mu_5$ and in the third using the library only contains the terms $\mu_1$ and $\mu_6$.}
\label{tab: cfr bb}
\end{table}



Finally, for the sake of completeness we show more details in Figure \ref{fig: Viscous Burgers bb_full} on the results, obtained with a black box, where the whole library was used.  The top left panel shows the solution and the top right panel  shows the error on parameters. In the bottom left panel, we show a comparison between the control found with the black box and the control computed by the algorithm that uses the implicit formula. The last plot shows a comparison between the costs of the controlled solution and the two RL solutions (implicit and black box). Note that, even if the model parameters found with the black box algorithm are less accurate than the ones found with the implicit one, the cost of the applied control is very similar.

\subsubsection{Test 3: Korteweg-De Vries}
In the third model we study the well-known Korteweg-De Vries (KdV) equation, with an additional diffusion term, which reads:
\begin{equation*}
\left\{ \begin{aligned}
    y_t(t,\xi) &= \frac{1}{2}y_{\xi\xi}(t,\xi) + 6 y(t,\xi)y_\xi(t,\xi) - y_{\xi\xi\xi}(t,\xi) + \chi_{[1,4]}(\xi) u(t), \qquad &&t \in [0,2], \ \xi \in (-10,7),\\
    y(0,\xi) &=\chi_{[0,6]}(\xi)\biggl(\cos\Bigl(\frac{\pi}{3}(\xi-3)\Bigr)+1\biggr)\qquad &&\xi \in 
    (-10,7),\\
    y(t,-10) &= 0, \ y(t,7) = 0, \qquad &&t \in [0,2].
    \end{aligned} \right.
\end{equation*}
Thus, it is a special case of \eqref{eq: parametric PDE} when $\mu_1^* = 0.5$, $\mu_6^* = 6$, $\mu_7^* = -1$. Note that in this test there is a third derivative in the equation. The boundary conditions are of Dirichlet type and the relative noise added was 1\%. The finite difference discretization is performed choosing $\Delta \xi = 0.1$ which leads a problem \eqref{eq: dynamics} of dimension $d=171$, and integrated in time with $\Delta t = 0.025.$
The configuration found from our Algorithm \ref{alg:RL} can be seen in Table \ref{tab: Kdv}. The update of the estimated configuration stopped after t=1.275 (i.e. after 51 iterations out of 80). The reconstructed parameter configuration has a difference of order $O(\Delta t)$ with respect to the true configuration. Again, as seen in the previous examples, this confirms the accurateness of our method. 

\begin{table}[thb]
\centering
\begin{tabular}{r|r|r|r|r|r|r|r}
True $\mu^*$ & 0.5 & 0 & 0 & 0 & 0 & 6 & -1\\ 
\hline
Estimated $\wt\mu$ & 0.4931 & 0.0012 & 0.0004 & 0.001 & -0.0016 & 5.9943 & -0.9999
\end{tabular}
\caption{Test 3: Korteweg-De Vries, $\Delta \xi = 0.1$, $\Delta t = 0.025$, 1\% relative noise.} 
\label{tab: Kdv}
\end{table}



\begin{figure}[htp]
\centering
\includegraphics[width=.3\textwidth]{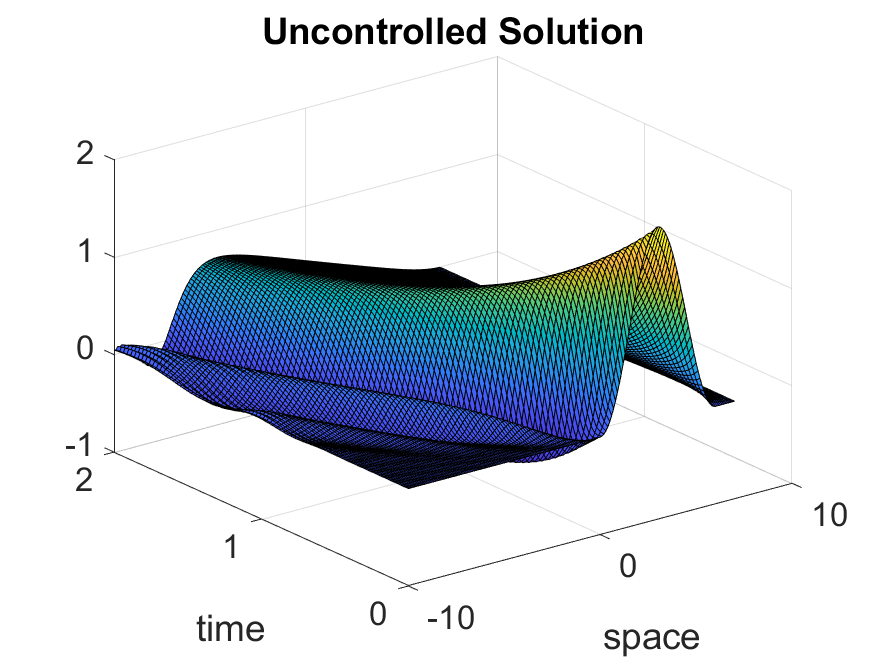}\hfill
\includegraphics[width=.3\textwidth]{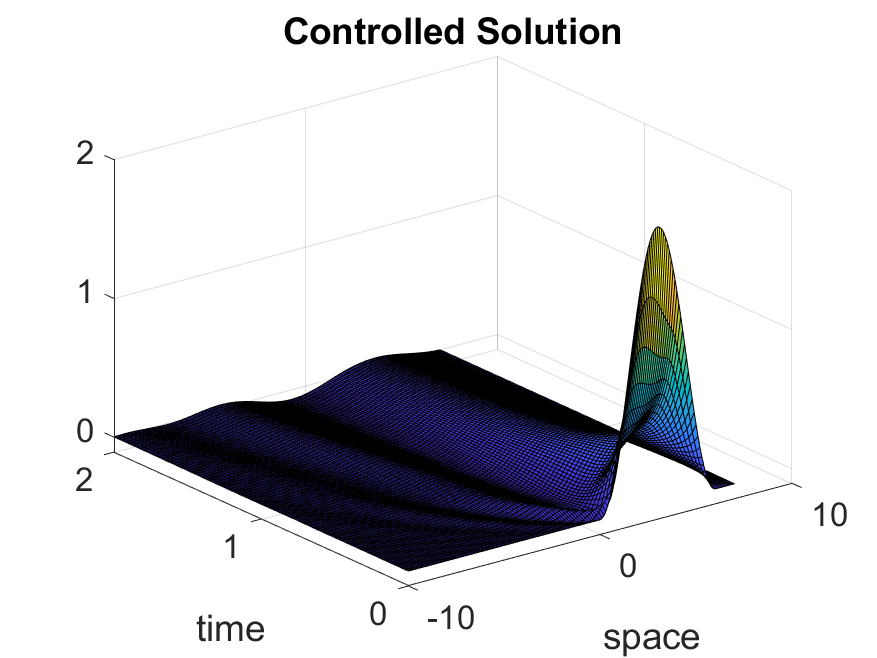}\hfill
\includegraphics[width=.3\textwidth]{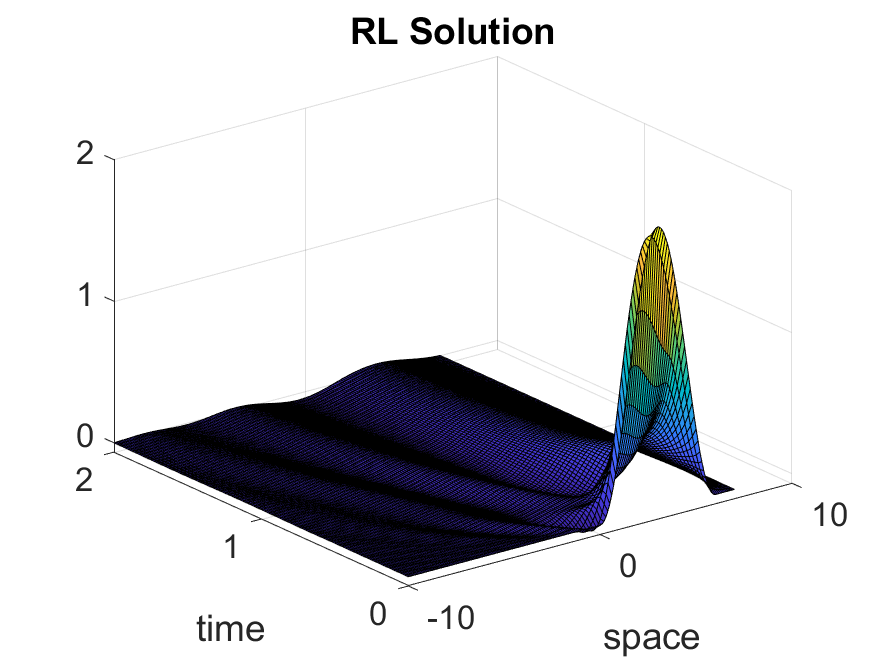}
\caption{Test 3: Korteweg-De Vries, $\Delta \xi = 0.1$, $\Delta t = 0.025$, 1\% relative noise.}
\label{fig: Kdv}
\end{figure}

The trajectories are presented in Figure \ref{fig: Kdv}. One can see that the middle and right panels have a similar behavior whereas the uncontrolled simulation has a completely different evolution.

\begin{figure}[htp]
\centering
\includegraphics[width=.3\textwidth]{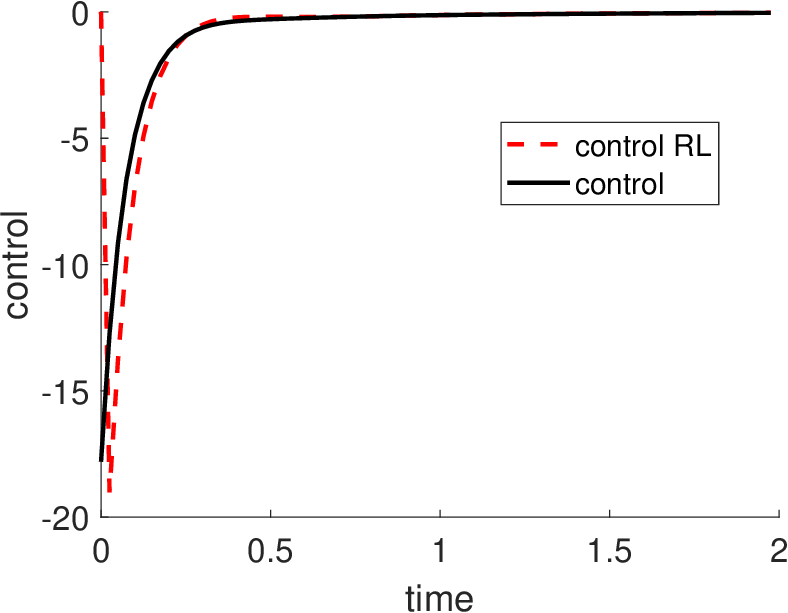}\hfill
\includegraphics[width=.3\textwidth]{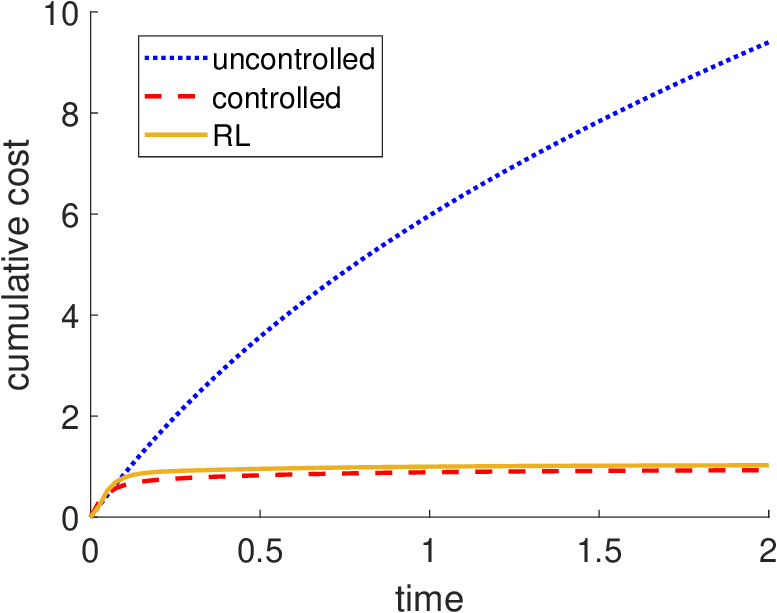}\hfill
\includegraphics[width=.3\textwidth]{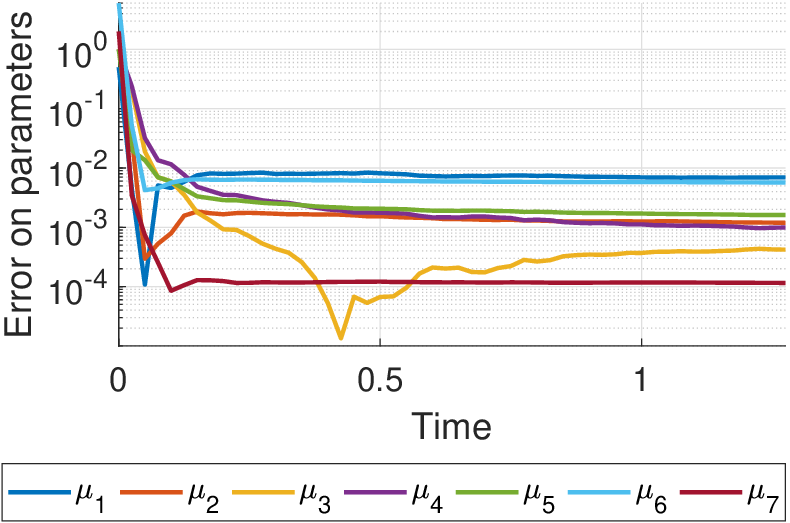}
\caption{Test 3: Korteweg-De Vries, $\Delta \xi = 0.1$, $\Delta t = 0.025$, 1\% relative noise. On the left the comparison between the control found using knowledge of the true $\mu$ and the control found by the RL algorithm is shown. In the middle the cumulative cost. On the right, the error on the parameter estimation at each time until the update stop.}
\label{fig: Kdv u mu}
\end{figure}

Finally, we show the computed controls in the left panel of Figure \ref{fig: Kdv u mu} which, after the first iterations, follow the same behavior. A more qualitative result is given in the middle panel of Figure \ref{fig: Kdv u mu} where we can see the evaluation of cost functional. Again (and as expected) Algorithm \ref{alg: sdre} performs slightly better than \ref{alg:RL} but still very close. To finalize, the history of the parameter configuration is shown in the right panel.

\subsection{CPU times}

In this subsection, we report in Table \ref{tab: CPU time}  the CPU times of the tests presented above. We compare the time needed to compute the uncontrolled, controlled and RL solutions for each of the three presented test cases for the implicit scheme. 

Since we have random components, the table has been obtained executing the algorithm 50 times and then considering the arithmetic mean of the execution times. Table \ref{tab: CPU time} shows that the time needed to obtain the solution with Algorithm \ref{alg: sdre} is similar to the time needed with our proposed method. This is because the number of PDEs solved is the same, the computation of the bayesian linear regression is neglectable since we do not deal with large scale problems and in our problem we have to solve one ARE less than SDRE since at the first iteration we decide to start with $0$ control. In the third test our method is slightly faster than SDRE, this also depends on the number of iteration needed in the Newton method which may be different since we opt for different control strategies.

\begin{table}[htbp]
    \centering
        \begin{tabular}{c|c|c|c}
        & uncontrolled & Algorithm \ref{alg: sdre} & Algorithm \ref{alg:RL}\\ 
        \hline
        Test 1 & $0.69$s & $8.7$s & $10.1$s\\ 
        \hline
        Test 2 & $0.87$s & $19.9$s & $21.3$s\\
        \hline
        Test 3 & $12.1$s & $67.7$s & $64.2$s\\ 
    \end{tabular}
    \caption{CPU times in seconds of the three presented tests. The times have been computed as the arithmetic mean of the time required to complete 50 algorithm's executions.}
\label{tab: CPU time}
\end{table}

To make the comparison fair, we consider the time needed to approximate the PDE in each method. Theoretically, one could think the black box in our method as an offline strategy with no cost. 



\begin{figure}[htbp]
\centering
\includegraphics[width=.3\textwidth]{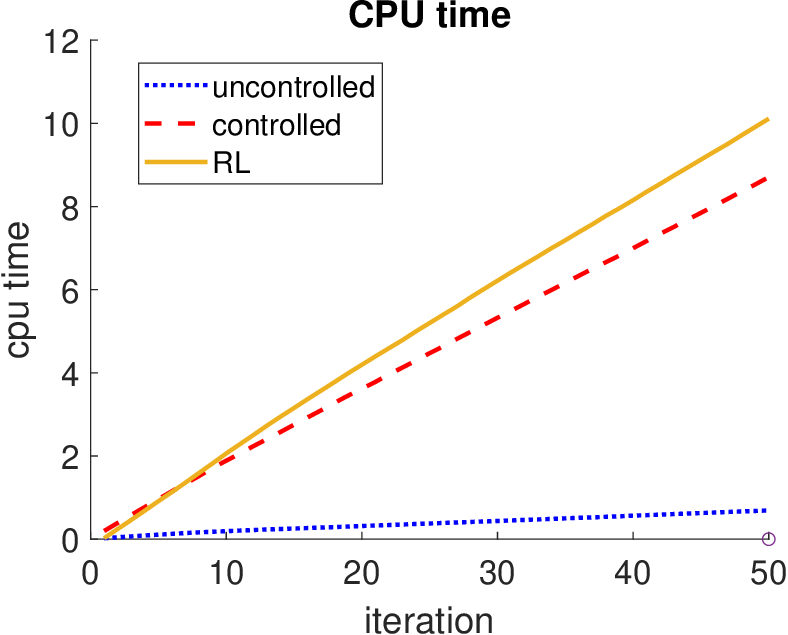}\hfill
\includegraphics[width=.3\textwidth]{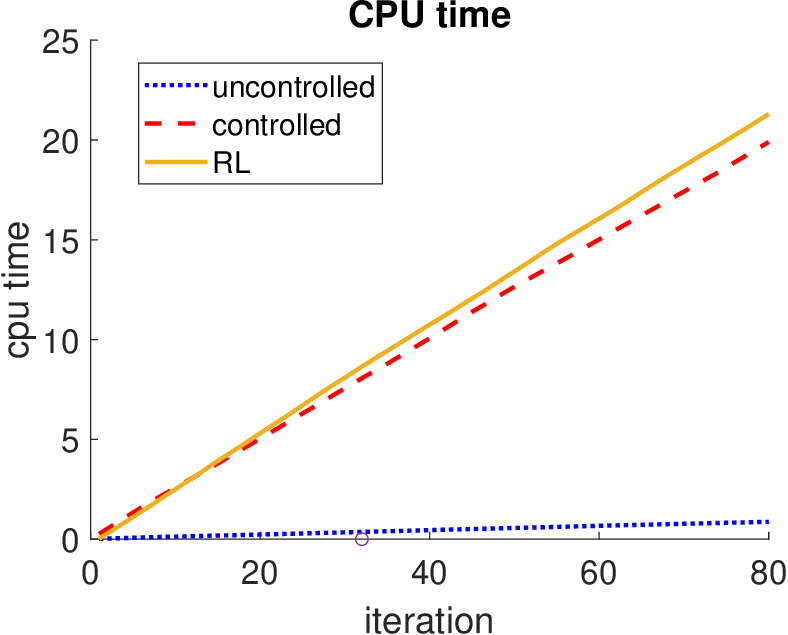}\hfill
\includegraphics[width=.3\textwidth]{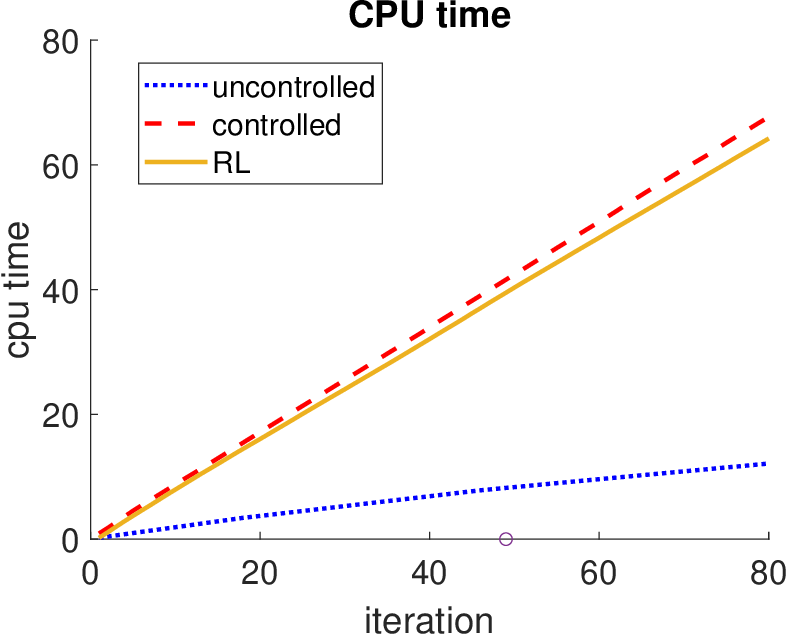}
\caption{Cumulative execution times at each iteration for Test 1 (left), Test 2 (middle) and Test 3 (Right). Mean times over 50 executions are considered.}
\label{fig: CPU time}
\end{figure}

Figure \ref{fig: CPU time} shows the execution time needed to conclude each iteration of the solution (and control) computation for the uncontrolled, controlled and RL cases. The final iteration times correspond to the times in Table \ref{tab: CPU time}. The uncontrolled case only requires the solution computation. We can observe that, in the first two tests, at the beginning the RL algorithm is slightly faster than the controlled one, and this is due to the choice of using a fixed control at the first iteration. Then, RL algorithm iterations are slightly slower, since more operations are carried out (e.g. bayesian regression). This behaviour is different for Test 3 as already commented. We also note that the uncontrolled KdV problem takes more time than the other two examples. 


\subsection{Convergence to the PDE}

\begin{figure}[htp]
\centering

\includegraphics[width=.31\textwidth]{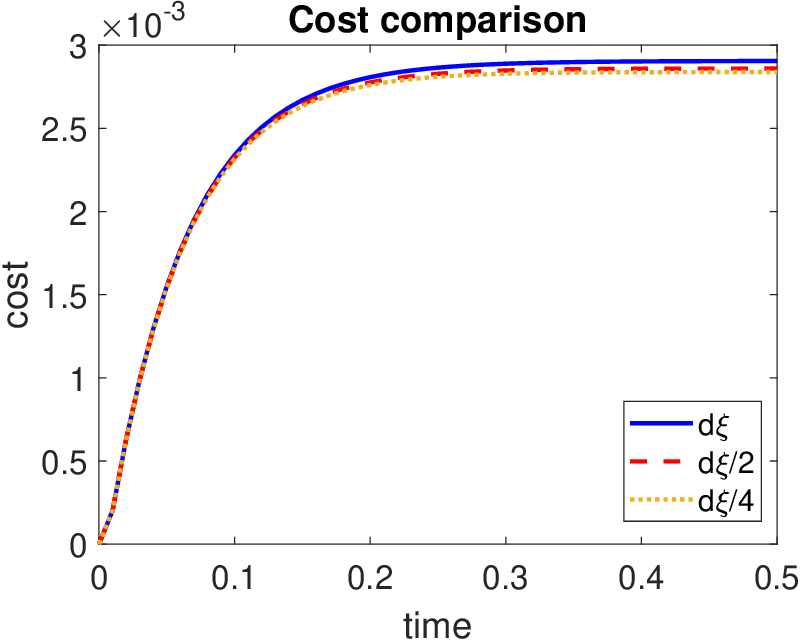}
\includegraphics[width=.31\textwidth]{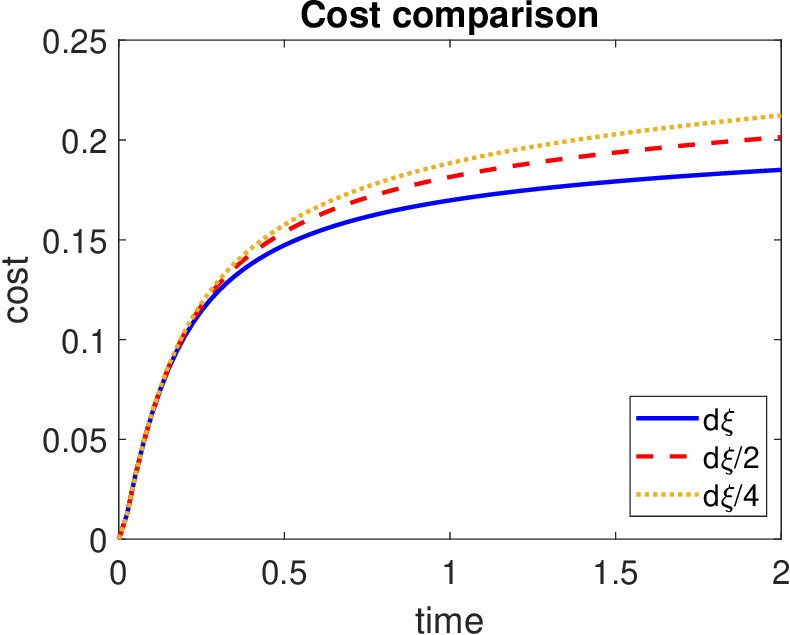}
\includegraphics[width=.31\textwidth]{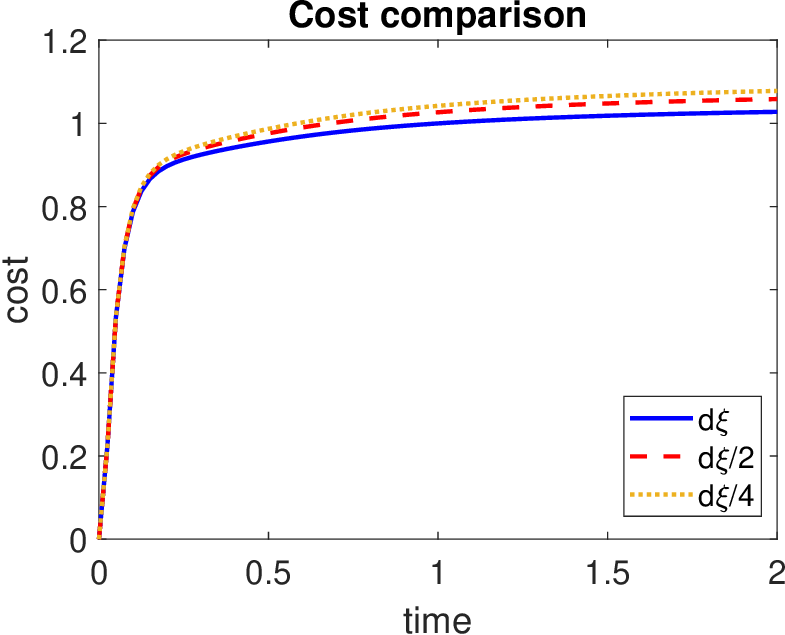}
\caption{Comparison of the cost functionals for Test 1 (left), Test 2(middle), Test 3(right)}
\label{fig:costdx}
\end{figure}

\begin{figure}[htp]
\centering
\includegraphics[width=.31\textwidth]{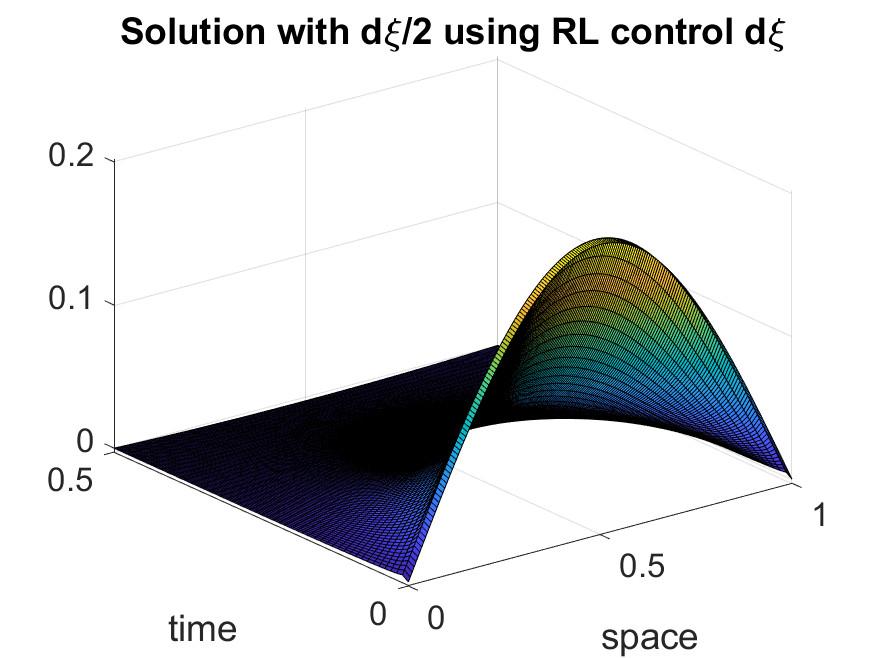}
\includegraphics[width=.31\textwidth]{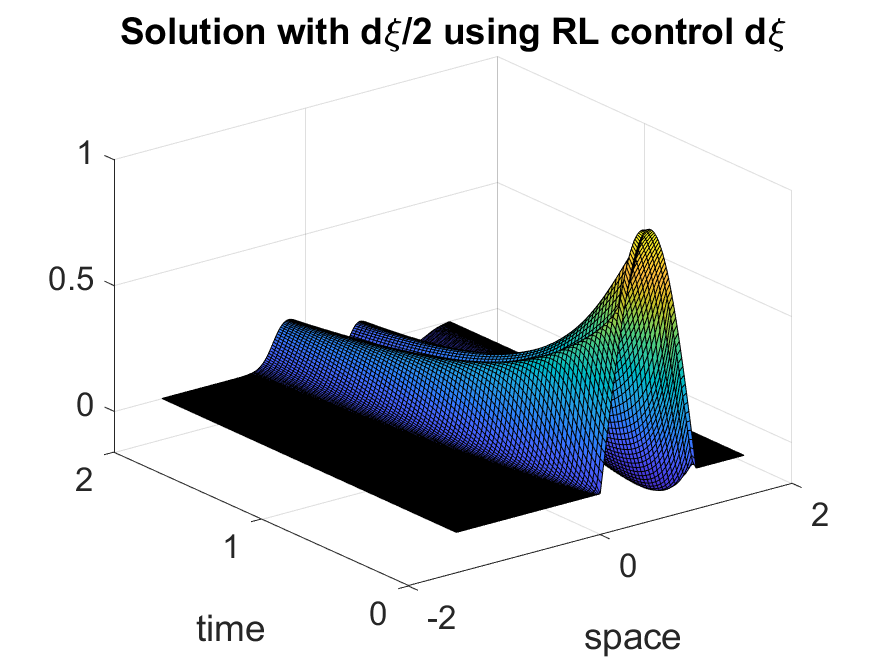}
\includegraphics[width=.31\textwidth]{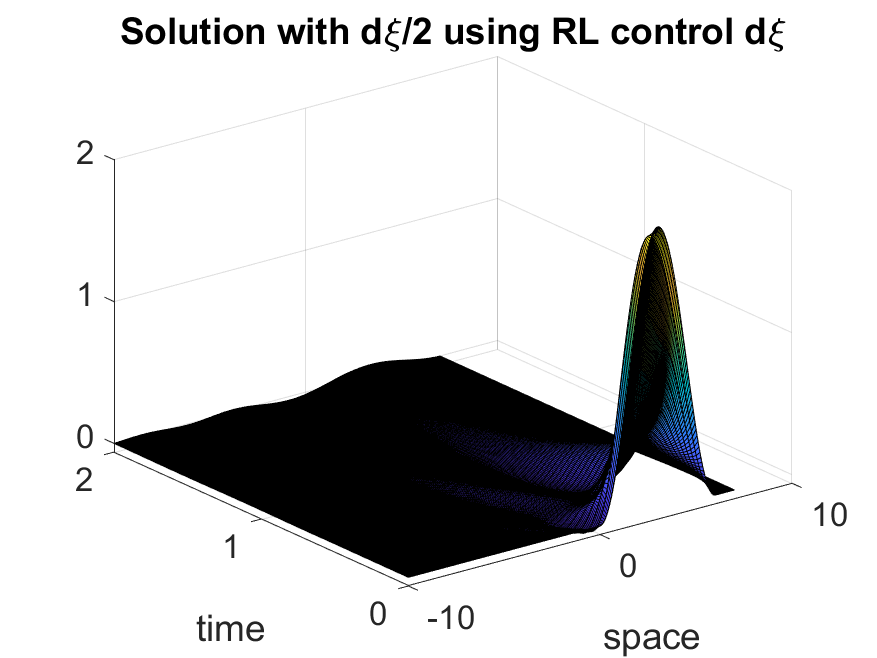}
\caption{Trajectories for Test 1 (left), Test 2(middle), Test 3(right)}
\label{fig:trajdx}
\end{figure}

\begin{figure}[htp]
\centering
\includegraphics[width=.31\textwidth]{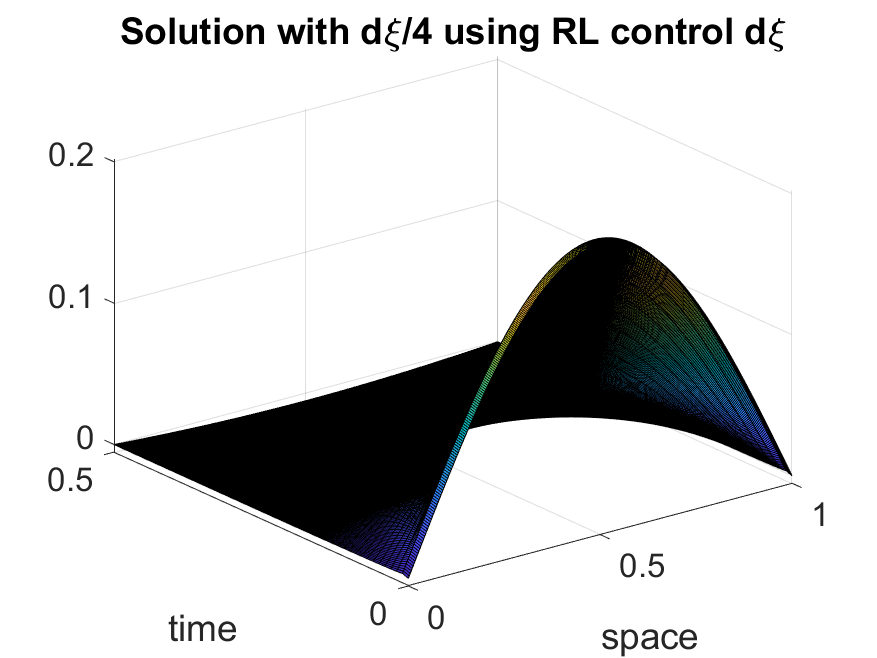}
\includegraphics[width=.31\textwidth]{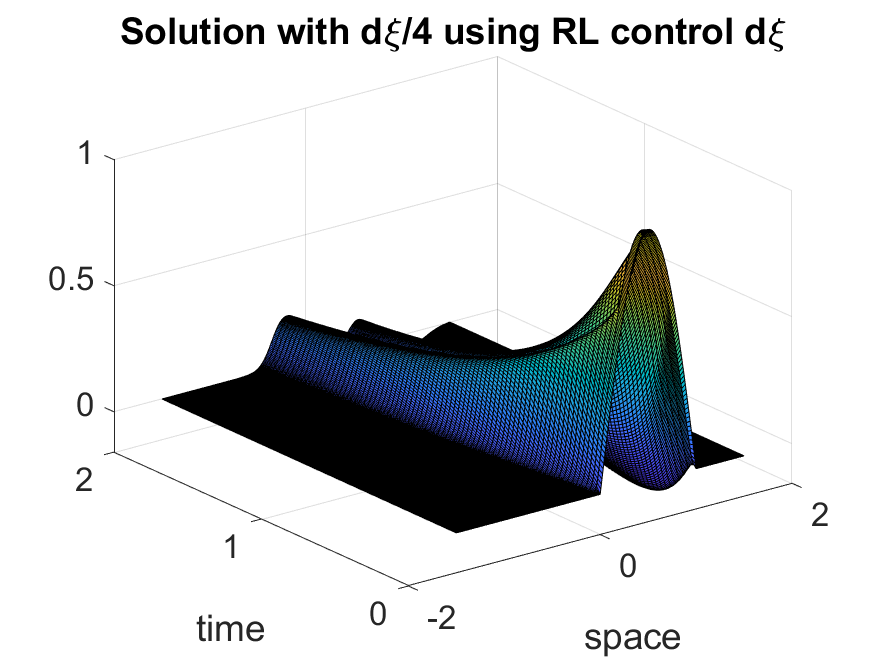}
\includegraphics[width=.31\textwidth]{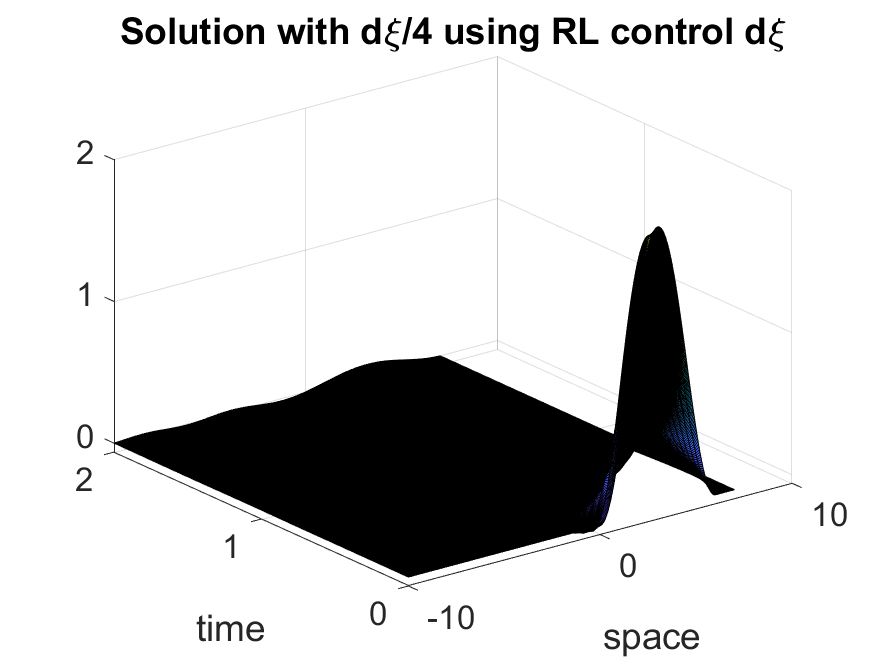}
\caption{Trajectories for Test 1 (left), Test 2(middle), Test 3(right)}
\label{fig:trajdx4}
\end{figure}

To conclude, we provide a numerical assessment of the convergence of our method in a PDE control framework. We consider the examples of the previous sections and study the convergence of the control for increasing dimension of the problem using the same time discretization grid used for each example to  study the role of the mesh towards the control of the PDE. Thus, we have tested the control obtained for a discretized problem of dimension $d$ (step $\Delta \xi$) using our Algorithm~\ref{alg:RL} and plugged into finer discretizations of the reference PDE of dimension, say $2d$ (step $\Delta \xi/2$) and $4d$ (step $\Delta \xi/4$). 
This has been done because, even if we use the the true parameter $\mu^*$ for the evolution, the obtained dynamics is still an approximation of the PDE evolution, due to the use of a numerical schemes. Finer grids allow us to better investigate the behaviour of the system after the application of the computed control.
For all the three numerical examples we plot the 3D solution generated with the finer grids (Figures \ref{fig:trajdx} and \ref{fig:trajdx4}) and the cost computed accordingly (Figure~\ref{fig:costdx}). We can see that the control found stabilizes the system also in these cases. 

\section{Conclusions}
\label{sec:conclusions}

We proposed a new algorithm designed to control/stabilize unknown PDEs under certain assumptions. The strength of the method is the identification of the system on the fly, where at each iteration we provide parameter estimate of the unknown system by Bayesian Linear regression. The update of the parameter configuration is based on the RL assumption where the user is always able to observe the true system configuration without its explicit knowledge.
Numerical experiments have shown convergent results that validate our proposed approach. Since, to the best of the authors' knowledge, this is the first approach of this kind for nonlinear problems, we leave several open problems, such as efficient algorithms for higher dimensional problems and a theoretical study of the convergence for the proposed method.
Then, it will be interesting to add further unknowns in the problem, such as e.g. the $B(x)$ term in the model and the quantity $Q$ in the cost.


 
\section*{Declarations}

\paragraph{\bf Funding} A. Alla and A. Pacifico are members of the INdAM-GNCS activity group. The work of M. Palladino is partially funded by the University of L'Aquila Starting Project Grant ``Optimal Control and Applications", and by INdAM-GNAMPA project, n. CUP\_E53C22001930001.

\paragraph{\bf Data Availibility} Enquires about code availability should be directed to the authors. No data has been used in this paper.

\paragraph{\bf Conflict of interest} No conflict of interest.

\paragraph{\bf Acknowledgments}
The authors want to express their deep gratitude to Maurizio Falcone. Thanks to him the authors met up and started to collaborate on this project. 

\bibliographystyle{plain}
\bibliography{bibsdre,Pesare_PhD_Thesis}

\begin{thebibliography}{10}

\bibitem{AKS23}
A.~Alla, D.~Kalise, and V.~Simoncini.
\newblock State-dependent riccati equation feedback stabilization for nonlinear
  pdes.
\newblock {\em Advances in Computational Mathematics}, 49, 2023.

\bibitem{AGW10}
Nils Altm\"{u}ller, Lars Gr\"{u}ne, and Karl Worthmann.
\newblock Receding horizon optimal control for the wave equation.
\newblock In {\em 49th IEEE Conference on Decision and Control (CDC)}, pages
  3427--3432, 2010.

\bibitem{BLT07}
H.~T. Banks, B.~M. Lewis, and H.~T. Tran.
\newblock Nonlinear feedback controllers and compensators: a state-dependent
  riccati equation approach.
\newblock {\em Comput. Optim. Appl.}, 37(2):177--218, Jun 2007.

\bibitem{Bardi-Capuzzo}
M.~Bardi and I.~Capuzzo-Dolcetta.
\newblock {\em Optimal Control and Viscosity Solutions of
  Hamilton-Jacobi-Bellman Equations}.
\newblock Birkhauser, 1997.

\bibitem{bellman1954}
Richard Bellman.
\newblock The theory of dynamic programming.
\newblock {\em Bulletin of the American Mathematical Society}, 60(6):503--515,
  1954.

\bibitem{B61}
Richard {B}ellman.
\newblock {\em Adaptive control processes: {A} guided tour}.
\newblock Princeton University Press, Princeton, N.J., 1961.

\bibitem{BH18}
P.~Benner and J.~Heiland.
\newblock Exponential stability and stabilization of extended linearizations
  via continuous updates of riccati-based feedback.
\newblock {\em International Journal of Robust and Nonlinear Control},
  28(4):1218--1232, 2018.

\bibitem{bertsekas2019reinforcement}
Dimitri Bertsekas.
\newblock {\em Reinforcement and Optimal Control}.
\newblock Athena Scientific, 2019.

\bibitem{bertsekas2008approximate}
Dimitri~P Bertsekas.
\newblock Approximate dynamic programming.
\newblock 2008.

\bibitem{box2011bayesian}
George~EP Box and George~C Tiao.
\newblock {\em Bayesian inference in statistical analysis}, volume~40.
\newblock John Wiley \& Sons, 2011.

\bibitem{BPK16}
S.L. Brunton, J.L. Proctor, and J.N. Kutz.
\newblock Discovering governing equations from data by sparse identification of
  nonlinear dynamical systems.
\newblock {\em Proceedings of the National Academy of Sciences of the United
  States of America,}, 115:3932–3937, 2016.

\bibitem{CFDKMRV21}
Sabrina Casper, Doris~H. Fuertinger, Peter Kotanko, Luca Mechelli, Jan Rohleff,
  and Stefan Volkwein.
\newblock Data-driven modeling and control of complex dynamical systems arising
  in renal anemia therapy.
\newblock In Matthias Ehrhardt and Michael G{\"u}nther, editors, {\em Progress
  in Industrial Mathematics at ECMI 2021}, pages 155--161, Cham, 2022. Springer
  International Publishing.

\bibitem{C97}
J.~R. {Cloutier}.
\newblock State-dependent {R}iccati equation techniques: an overview.
\newblock In {\em Proceedings of the 1997 American Control Conference (Cat.
  No.97CH36041)}, volume~2, pages 932--936 vol.2, 1997.

\bibitem{falcone2013semi}
Maurizio Falcone and Roberto Ferretti.
\newblock {\em Semi-Lagrangian approximation schemes for linear and
  Hamilton—Jacobi equations}.
\newblock SIAM, 2013.

\bibitem{freedman2009statistical}
David~A Freedman.
\newblock {\em Statistical models: theory and practice}.
\newblock cambridge university press, 2009.

\bibitem{larsnmpc}
Lars Gr\"{u}ne and J\"{u}rgen Pannek.
\newblock {\em Nonlinear model predictive control}.
\newblock Communications and Control Engineering Series. Springer, London,
  2011.
\newblock Theory and algorithms.

\bibitem{SAC2018}
Tuomas Haarnoja, Aurick Zhou, Pieter Abbeel, and Sergey Levine.
\newblock Soft actor-critic: Off-policy maximum entropy deep reinforcement
  learning with a stochastic actor.
\newblock In {\em International Conference on Machine Learning}, pages
  1861--1870. PMLR, 2018.

\bibitem{KKB19}
E.~Kaiser, J.~N. Kutz, and S.~L. Brunton.
\newblock Sparse identification of nonlinear dynamics for model predictive
  control in the low-data limit.
\newblock {\em Proceedings of the Royal Society A: Mathematical, Physical and
  Engineering Sciences}, 474(2219):20180335, 2018.

\bibitem{KKLPWY21}
G.E. Karniadakis, I.G. Kevrekidis, L.~Lu, P.~Perdikaris, S.~Wang, and L.~Yang.
\newblock Physics-informed machine learning.
\newblock {\em Nat Rev Phys}, 3:686--707, 2021.

\bibitem{KK04}
D.~A. Knoll and D.~E. Keyes.
\newblock Jacobian-free {N}ewton-{K}rylov methods: a survey of approaches and
  applications.
\newblock {\em J. Comput. Phys.}, 193(2):357--397, 2004.

\bibitem{KS08}
Miroslav Krstic and Andrey Smyshlyaev.
\newblock Adaptive control of pdes.
\newblock {\em IFAC Proceedings Volumes}, 40(13):20--31, 2007.
\newblock 9th IFAC Workshop on Adaptation and Learning in Control and Signal
  Processing.

\bibitem{L07}
Randall~J. LeVeque.
\newblock {\em Finite difference methods for ordinary and partial differential
  equations}.
\newblock Society for Industrial and Applied Mathematics (SIAM), Philadelphia,
  PA, 2007.
\newblock Steady-state and time-dependent problems.

\bibitem{Lillicrap2016}
Timothy~P Lillicrap, Jonathan~J Hunt, Alexander Pritzel, Nicolas Heess, Tom
  Erez, Yuval Tassa, David Silver, and Daan Wierstra.
\newblock Continuous control with deep reinforcement learning.
\newblock In {\em 4th International Conference on Learning Representations
  (ICLR)}, 2016.

\bibitem{MLG20}
Andreas~B. Martinsen, Anastasios~M. Lekkas, and Sébastien Gros.
\newblock Combining system identification with reinforcement learning-based
  mpc.
\newblock {\em IFAC-PapersOnLine}, 53(2):8130--8135, 2020.
\newblock 21st IFAC World Congress.

\bibitem{Mnih2015}
Volodymyr Mnih, Koray Kavukcuoglu, David Silver, Andrei~A Rusu, Joel Veness,
  Marc~G Bellemare, Alex Graves, Martin Riedmiller, Andreas~K Fidjeland, Georg
  Ostrovski, et~al.
\newblock Human-level control through deep reinforcement learning.
\newblock {\em Nature}, 518(7540):529--533, 2015.

\bibitem{pacifico2021LSSC}
Agnese Pacifico, Andrea Pesare, and Maurizio Falcone.
\newblock A new algorithm for the {LQR} problem with partially unknown
  dynamics.
\newblock In {\em 2021 Large-Scale Scientific Computing (LSSC)}. Springer
  International Publishing, in press.

\bibitem{powell2007approximate}
Warren~B Powell.
\newblock {\em Approximate Dynamic Programming: Solving the curses of
  dimensionality}, volume 703.
\newblock John Wiley \& Sons, 2007.

\bibitem{powell2021reinforcement}
Warren~B Powell.
\newblock From reinforcement learning to optimal control: A unified framework
  for sequential decisions.
\newblock In {\em Handbook of Reinforcement Learning and Control}, pages
  29--74. Springer, 2021.

\bibitem{RPK19}
M.~Raissi, P.~Perdikaris, and G.E. Karniadakis.
\newblock Physics-informed neural networks: A deep learning framework for
  solving forward and inverse problems involving nonlinear partial differential
  equations.
\newblock {\em Journal of Computational Physics}, 378:686--707, 2019.

\bibitem{Rasmussen2006}
CE. Rasmussen and CKI. Williams.
\newblock {\em Gaussian Processes for Machine Learning}.
\newblock Adaptive Computation and Machine Learning. MIT Press, Cambridge, MA,
  USA, January 2006.

\bibitem{rossi2012bayesian}
Peter~E Rossi, Greg~M Allenby, and Rob McCulloch.
\newblock {\em Bayesian statistics and marketing}.
\newblock John Wiley \& Sons, 2012.

\bibitem{RABK19}
Samuel Rudy, Alessandro Alla, Steven~L. Brunton, and J.~Nathan Kutz.
\newblock Data-driven identification of parametric partial differential
  equations.
\newblock {\em SIAM Journal on Applied Dynamical Systems}, 18(2):643--660,
  2019.

\bibitem{RBPK16}
S.H. Rudy, S.L. Brunton, J.L. Proctor, and J.N. Kutz.
\newblock Data-driven discovery of partial differential equations.
\newblock {\em Science Advances}, 3, 2017.

\bibitem{rummery1994line}
Gavin~A Rummery and Mahesan Niranjan.
\newblock {\em On-line Q-learning using connectionist systems}, volume~37.
\newblock Citeseer, 1994.

\bibitem{TRPO2015}
John Schulman, Sergey Levine, Pieter Abbeel, Michael Jordan, and Philipp
  Moritz.
\newblock Trust region policy optimization.
\newblock In {\em International conference on machine learning}, pages
  1889--1897. PMLR, 2015.

\bibitem{sutton1988learning}
Richard~S Sutton.
\newblock Learning to predict by the methods of temporal differences.
\newblock {\em Machine learning}, 3(1):9--44, 1988.

\bibitem{SuttonBarto1edn}
Richard~S Sutton and Andrew~G Barto.
\newblock {\em Reinforcement learning: An introduction}, volume~1.
\newblock MIT Press, Cambridge, MA, first edition, 1998.

\bibitem{SuttonBarto}
Richard~S Sutton and Andrew~G Barto.
\newblock {\em Reinforcement Learning: An Introduction}.
\newblock MIT Press, Cambridge, MA, second edition, 2018.

\bibitem{watkins1989learning}
Christopher Watkins and John~Cornish Hellaby.
\newblock Learning from delayed rewards.
\newblock 1989.

\end{thebibliography}



\end{document}